\documentclass[12pt,draftcls,onecolumn]{IEEEtran}

\IEEEoverridecommandlockouts                              

\makeatletter
\let\NAT@parse\undefined
\makeatother

\usepackage{graphicx}
\usepackage{bm}
\usepackage{epsfig}
\usepackage{amsthm}
\usepackage{amssymb,amsfonts,amsmath}
\usepackage{subfigure}
\usepackage{setspace}
\usepackage{enumerate}
\usepackage[numbers]{natbib}
\usepackage{color}
\usepackage{dsfont}
\usepackage{verbatim}

\newtheorem{theorem}{Theorem}
\newtheorem{lemma}{Lemma}
\newtheorem{remark}{Remark}
\newtheorem{defn}{Definition}
\newtheorem{proposition}{Proposition}
\newtheorem{example}{Example}

\newtheorem{corollary}{Corollary}

\def\<{\langle}
\def\>{\rangle}

\def \b {\beta}

\def \l {\lambda}
\def \e {\varepsilon}

\def \F {\mathcal{F}}

\title{\LARGE \bf
On Separating Points for Ensemble Controllability
}

\author{Jr-Shin Li, Wei Zhang, and Lin Tie
\thanks{*This work was supported in part by the National Science Foundation under the awards CMMI-1462796, ECCS-1509342, and ECCS-1810202, and by the Air Force Office of Scientific Research under the award FA9550-17-1-0166.}
\thanks{J.-S. Li is with the Department of Electrical and Systems Engineering, Washington University,
        St. Louis, MO 63130, USA
        {\tt\small jsli@wustl.edu}}%
\thanks{W. Zhang is with the Department of Electrical and Systems Engineering, Washington University,
		St. Louis, MO 63130, USA
        {\tt\small wei.zhang@wustl.edu}}%
\thanks{L. Tie is with the School of Automation Science and Electrical Engineering, Beihang University, Beijing, 100191, China
        {\tt\small tielin@buaa.edu.cn}}%
}

\begin{document}

\maketitle




\begin{abstract}
Recent years have witnessed a wave of research activities in systems science toward the study of population systems. The driving force behind this shift was geared by numerous emerging and ever-changing technologies in life and physical sciences and engineering, from neuroscience, biology, and quantum physics to robotics, where many control-enabled applications involve manipulating a large ensemble of structurally identical dynamic units, or agents. Analyzing fundamental properties of ensemble control systems in turn plays a foundational and critical role in enabling and, further, advancing these applications, and the analysis is largely beyond the capability of classical control techniques. In this paper, we consider an ensemble of time-invariant linear systems evolving on an infinite-dimensional space of continuous functions. We exploit the notion of separating points and techniques of polynomial approximation to develop necessary and sufficient ensemble controllability conditions. In particular, we introduce an extended notion of controllability matrix, called \emph{Ensemble Controllability Gramian}. This means enables the characterization of ensemble controllability through evaluating controllability of each individual system in the ensemble. As a result, the work provides a unified framework with a systematic procedure for analyzing control systems defined on an infinite-dimensional space by a finite-dimensional approach.
\end{abstract}




\section{Introduction}
\label{sec:intro}
The emergence of complex systems constituted by an ensemble of structurally identical dynamical units (agents) has made waves driving the recent research in systems science towards studying the control of collective dynamics and behavior of population systems \cite{Rosenblum04,kiss02,YX_Chen_thesis,YX_Chen16,Politi2015}. Notable examples include stimulating neuronal populations to alleviate the symptoms of neurological disorders such as Parkinson’s disease \cite{Ching2013b,Li_TAC13,Kafashan2015}, using optogenetics to manipulate circadian rhythms in the suprachiasmatic nucleus (SCN) network \cite{Sidor14}, controlling synchronization patterns in a network of nonlinear oscillators \cite{Rosenblum04,Li_NatureComm16,Kiss:2007:1886-1889}, steering a collection of robots under model perturbation \cite{Becker12}, and exciting spin ensembles for spectroscopy, imaging, and quantum computation \cite{Glaser98,Li_PRA06,Li_PNAS11}. The fundamental difficulty in these applications involving dynamic populations is that control and observation can only be implemented at the population level, i.e., through broadcasting a single input signal to all the systems in the population and through receiving aggregated measurements of the systems in the population, respectively. This restriction gives rise to a new control paradigm of population-based control, called ensemble control. The principles and technologies of ensemble control remain far from reaching their full potential and capacity.

Thanks to the richness on the subject of ensemble control, in the past decade research in ensemble systems took a sharp turn, and a great deal of focus was placed on investigating fundamental properties of ensemble systems and devising tractable computational algorithms for calculating optimal ensemble control laws. Various specialized methods that are nontrivially linked to control theory, such as polynomial approximation, complex functional analysis, and statistically moment-based approaches, were developed to study ensemble controllability of linear \cite{Li_TAC11,Schonlein13,Helmke14,Li_TAC16,Zeng_16_moment,Dirr16,Li_SCL18}, bilinear \cite{Li_TAC09,Belhadj15}, and some forms of nonlinear \cite{Li_TAC13} ensemble systems. Algebraic geometry and probabilistic methods were established to formulate and analyze problems of ensemble observability, where the innovation was to map the observability analysis to a problem in a remote domain, such as mathematical tomography or representation theory \cite{zeng2016tac,zeng2015nonlinear_ensemble,XD_Chen19}. In addition, the complexity and tremendous scale of ensemble systems has pushed the boundaries of computational optimal control toward developing customized and effective numerical algorithms for solving optimal control problems involving ensemble systems, such as iterative methods \cite{Li_SICON17,Li_Automatica17}, pseudospectral methods \cite{Li_TAC12_QCP,Phelps14}, operator-theoretic methods \cite{Zeng18}, and sample average approximations \cite{Gong16}. Moreover, these novel developments not only advanced our fundamental understanding of ensemble systems but also provided immediate control-enabled applications in highly diverse areas, including neuroscience, biology, quantum control, and complex networks \cite{Ching2013b,Kafashan2015,Li_NatureComm16,kuritz2019ensemble,kiss02,Li_NatureComm17,Li_PNAS11,Li_JCP11,Chen14,Jiang12,Augier18,Rosenblum04,Kiss:2007:1886-1889,Li_PNAS18}. 

Although much progress and viable analysis has been conducted on controllability of linear ensemble systems providing variant necessary and/or sufficient conditions, a transparent and unified view of what makes such an ensemble system controllable is missing. In this paper, we fill the gap and offer an algebraic explanation by exploiting the notion of \emph{separating points} originated from the theory of polynomial approximation \cite{Folland99}, which gives rise to a tangible geometric understanding of controllability of linear, and possibly broader classes of, ensemble systems. Based on this algebraic-geometric point of view, we develop verifiable necessary and sufficient ensemble controllability conditions for the general time-invariant linear ensemble system, indexed by a scalar parameter varying on some compact set, by using integrated techniques of separating points, polynomial approximation, and Lie algebra. In particular, we illuminate the effect of the functional property and the spectrum structure of the drift dynamics on ensemble controllability and establish a reparameterization approach to enable the analysis based on the idea of separating points. Moreover, we introduce an extended notion of controllability matrix, called \emph{Ensemble Controllability Gramian}, which enables a counterintuitive characterization of ensemble controllability through evaluating controllability of each individual system in the ensemble under a different parameterization. This work provides a unified framework with a systematic procedure for analyzing controllability of linear ensemble systems and control systems defined on an infinite-dimensional space using a finite-dimensional approach.

This paper is organized as follows. In Section \ref{sec:1D}, we introduce the notion of ensemble controllability and utilize techniques of separating points and polynomial approximation to construct uniform ensemble controllability conditions for an ensemble of scalar time-invariant linear systems. In Section \ref{sec:n-d}, we extend the results developed in Section \ref{sec:1D} to analyze uniform ensemble controllability of multi-dimensional linear ensemble systems. The main idea is to map such high-dimensional ensembles to a scalar ensemble system by using a spectrum-preserving coordinate transformation. We analyze ensemble system whose system matrix is diagonalizable or similar to a Jordan canonical form with real eigenvalues.

\section{Uniform ensemble controllability of one-dimensional linear ensemble systems}
\label{sec:1D}
In this section, we establish explicit uniform ensemble controllability conditions for the one-dimensional time-invariant linear ensemble system. In particular, we introduce the notion of ``separating points'' (see Definition \ref{def:separating_points} in Appendix \ref{appd:S-W_thm}) from the theory of polynomial approximation \cite{Apostol74} for analyzing this fundamental property. Through the analysis, we reveal the important role that the injectivity of the system dynamics and the number of independent control functions play for characterization of uniform ensemble controllability. 

\begin{defn}[Ensemble Controllability]
	\label{def:controllability}
	An ensemble of systems defined on a manifold $M$ parameterized by a parameter $\b$ taking values on a space $K$, given by 
	\begin{equation}
		\label{eq:ensemble}
		\frac{d}{dt}x(t,\b)=F(t,\b,x(t,\b),u(t)),
	\end{equation}
	where $x(t,\cdot)\in\mathcal{F}(K)$ is the state and $\mathcal{F}(K)$ is a space of $M$-valued functions defined on $K$, is said to be ensemble controllable on $\mathcal{F}(K)$, if for any $\e>0$ and starting with any initial state $x_0\in\mathcal{F}(K)$, where $x_0(\b)=x(0,\b)$, there exists a piecewise constant control signal $u:[0,T]\to\mathbb{R}^m$ that steers the system into an $\e$-neighborhood of a desired target state $x_F\in\mathcal{F}(K)$ at a finite time $T>0$, i.e., $d(x(T,\cdot),x_F(\cdot))<\e$, where $d:\mathcal{F}(K)\times \mathcal{F}(K)\rightarrow\mathbb{R}$ is a metric on $\mathcal{F}(K)$. 
\end{defn}

In particular, if we equip the state space $\F(K)$ of the ensemble system in \eqref{eq:ensemble} with the uniform metric, i.e., $d(f,g)=\sup_{\b\in K}\rho(f,g)$ for any $f,g\in\F(K)$, where $\rho:M\times M\rightarrow\mathbb{R}$ is a metric on $M$, Definition \ref{def:controllability} is referred to as \emph{uniform ensemble controllability}. 

In this paper, we focus on the case in which $K$ is a compact subset of $\mathbb{R}$, $\F(K)=C(K,\mathbb{R}^n)$ is the space of continuous $\mathbb{R}^n$-valued functions, and $d$ is induced by the uniform norm, i.e., $\|f\|_{\infty}=\sup_{\b\in K}\|f(\b)\|$ for each $f\in\F$, where $\|\cdot\|$ is any norm on $\mathbb{R}^n$. Our starting point for controllability studies is a systematic analysis for the one-dimensional linear ensemble system, from single-input to multi-input scenarios.

\subsection{Single-input linear ensembles}
\label{sec:single_input}
In this section, we analyze the one-dimensional single-input linear ensemble system and illuminate the role of monotonicity of the natural (drift) dynamics plays for ensuring controllability.  

\begin{proposition}
    \label{prop:1D}
    Consider the ensemble of one-dimensional single-input linear systems, indexed by the parameter $\b$ varying on a compact set $K=[\b_1,\b_2]\subset\mathbb{R}$,
    \begin{align}
        \label{eq:1D}
        \frac{d}{dt}x(t,\b)=a(\beta)x(t,\b)+b(\b)u(t),
    \end{align}
    where $x(t,\cdot)\in C(K,\mathbb{R})$, $a,b\in C(K,\mathbb{R})$, and $u:[0,T]\rightarrow\mathbb{R}$ is piecewise constant. This system is uniformly ensemble controllable on $C(K,\mathbb{R})$ if and only if 
	\begin{enumerate}
		\item[{\rm (i)}] $a$ is injective, and 
		\item[{\rm (ii)}] $b$ is nowhere vanishing 
			on $K$.
	\end{enumerate}
\end{proposition}

{\it Proof:} The reachable set of this system is characterized by the Lie algebra generated by its drift and control vector fields, 
that is, $\mathcal{L}_0={\rm span}\{a^kb:\ k=0,1,2,\ldots,\}$. Thus, by Definition \ref{def:controllability}, the system in \eqref{eq:1D} is uniformly ensemble controllable on $C(K,\mathbb{R})$ if and only if $\overline{\mathcal{L}_0}=C(K,\mathbb{R})$, i.e., for any given $\e>0$ and $\xi\in C(K,\mathbb{R})$, there exists $\eta\in\mathcal{L}_0$ such that $\|\xi-\eta\|_{\infty}<\e$, where $\eta(\b)=\sum_{k=0}^{N}c_ka^k(\b)b(\b)$ for some $N>0$ ($N$ may depend on $\e$) and $c_0,\dots,c_N\in\mathbb{R}$. Here, $\overline{\mathcal{L}_0}$ denotes the closure of $\mathcal{L}_0$ under the uniform topology. Consequently, uniform ensemble controllability of the system in \eqref{eq:1D} is equivalent to whether such a polynomial approximation can be achieved, i.e., the existence of $N>0$ and $c_0,\dots,c_N\in\mathbb{R}$.

(Sufficiency): 
Observe that $\|\xi-\sum_{k=0}^Nc_ka^{k}b\|_{\infty}<\e$ is equivalent to 
\begin{align}
    \label{eq:beta_approx}
\big\|\hat{\xi}-\sum_{k=0}^Nc_ka^k\big\|_{\infty}<\e',
\end{align}
provided $b(\b)\neq0$ for all $\b\in K$, where $\hat{\xi}(\b)=\xi(\b)/b(\b)$ and $\e'=\e/\|b\|_{\infty}$. Moreover, $\xi,b\in C(K,\mathbb{R})$ implies $\hat{\xi}\in C(K,\mathbb{R})$.

Let $S=a(K)\subset\mathbb{R}$ be the range of $a$, then the injectivity and continuity of $a:K\rightarrow S$ implies that it has a unique continuous inverse function defined on $S$, i.e., $a^{-1}:S\rightarrow K$ such that for any $\b\in K$, there exists an unique $s\in S$ such that $a^{-1}(s)=\b$. Because $\hat{\xi}\circ a^{-1}$ is a continuous function on $S$, the Weierstrass approximation theorem \cite{Apostol74} guarantees the possibility of the approximation, $\big\|(\hat{\xi}\circ a^{-1})(s)-\sum_{k=0}^Nc_ks^k\big\|_{\infty}<\e'$, which leads to \eqref{eq:beta_approx}. Therefore, the system in \eqref{eq:1D} is uniformly ensemble controllable on $C(K,\mathbb{R})$.

(Necessity):  
We will show that if either the injectivity of $a$ or non-vanishing of $b$ is violated, then the system in \eqref{eq:1D} is not uniformly ensemble controllable on $C(K,\mathbb{R})$. First, suppose that there exists $\b_0\in K$ such that $b(\b_0)=0$, then the system indexed by $\b_0$, i.e., $\frac{d}{dt}x(t,\b_0)=a(\b_0)x(t,\beta_0)+b(\b_0)u=a(\b_0)x(t,\beta_0)$, is not controllable on $\mathbb{R}$, so that the whole ensemble is not controllable on $C(K,\mathbb{R})$. This implies that condition (ii) 
is necessary. Next, if $a$ is not injective, then there exist $\bar{\b}_1,\bar{\b}_2\in K$ such that $\bar{\b}_1\neq\bar{\b}_2$ and $a(\bar{\b}_1)=a(\bar{\b}_2)$. Because every function $f\in\mathcal{L}_0$ is in  the form $f(\b)=\sum_{k=0}^Nc_ka^{k}(\b)b(\b)$ for some $N\in\mathbb{N}$ and constants $c_k\in\mathbb{R}$, it follows that the equality $f(\bar{\b}_1)/b(\bar{\b}_1)=f(\bar{\b}_2)/b(\bar{\b}_2)$ must hold, and the functions in $C(K,\mathbb{R})$ violating this relation cannot be approximated by the elements in $\mathcal{L}_0$. Therefore, the system fails to be uniformly ensemble controllable. This concludes that $a(\b)$ must be injective. \hfill$\Box$

An immediate consequence of the above result tells that the ensemble system with linear parameter variation in the natural dynamics, given by $\frac{d}{dt}x(t,\b)=\beta x(t,\b)+u(t)$, is uniformly ensemble controllable. Also, when the conditions in Proposition \ref{prop:1D} are violated, the system will not be ensemble controllable as illustrated in the following example.

\begin{example}
    \label{ex:1D_single}
    \rm Consider the scalar single-input linear ensemble system
	\begin{equation}
		\label{eq:1D_single}
		\frac{d}{dt}x(t,\b)=\b^2x(t,\b)+u(t),
	\end{equation}
    where $x(t,\cdot)\in C(K,\mathbb{R})$, $\b\in K=[-1,1]$, and $u:[0,T]\to\mathbb{R}$. According to Proposition \ref{prop:1D}, this system is not ensemble controllable since $\b^2$ is not injective on $K$. 
\end{example}

\begin{remark}
	\rm
	The uncontrollability of the system in \eqref{eq:1D_single} is also evident through the fact that the Lie algebra $\mathcal{L}_0$ characterizing the reachable set of this system, given by
	$$\mathcal{L}_0=\text{span}\{1,\b^2,\b^4,\ldots\},$$
	is not dense in $C(K,\mathbb{R})$, since any ensemble state $\xi(\b)$ that is an odd function on $K$ cannot be reached. In other words, for any $\b\in K$, there is no function in $\mathcal{L}_0$ that is able to separate the two points $\b$ and $-\b$, namely, the systems indexed by $\b$ and $-\b$ are not distinguishable by using the available single input (see Definition \ref{def:separating_points} in Appendix \ref{appd:S-W_thm} for the idea of separating points).
\end{remark}

However, when multiple control inputs are available, the injectivity of the drift dynamics $a(\b)$ is no longer necessary, since extra controls can be used to separate, or distinguish, points of the same image in different injective branches. This is studied in the next Section.

\subsection{Multi-input linear ensembles}
\label{sec:multi-input}
Distinct from the single-input scenario, in the case in which the ensemble system receives multiple inputs, non-injectivity of the drift dynamics failing controllability can be compensated for by using additional control inputs. 
The exemplar presented below inspires the need of considering the interplay between the injectivitiy of the system dynamics and the number of independent control fields in order to fully characterize controllability of ensemble systems.

\begin{example}
    \label{ex:1D_multiple}
    \rm Consider the scalar multi-input linear ensemble system
    $$\frac{d}{dt}x(t,\b)=\b^2x(t,\b)+u_1(t)+\b u_2(t),$$
    where $x(t,\cdot)\in C(K,\mathbb{R})$, $\b\in K=[-1,1]$, and $u_1,u_2: [0,T]\to\mathbb{R}$. Although $a(\b)=\b^2$ is not injective on $K$, this system is uniformly ensemble controllable. In this case, the reachable set is characterized by
    $$\mathcal{L}_0=\text{span}\{1,\b,\b^2,\b^3,\ldots\},$$
    and hence its closure $\overline{\mathcal{L}_0}=C(K,\mathbb{R})$. This implies that any given state $\xi(\b)$ of the ensemble is reachable, i.e., for any $\e>0$, there exist $c_k$ and an integer $N>0$ 
	such that $\|\xi(\b)-\sum_{k=0}^{N} c_k\b^k\|_{\infty}<\e$. 
\end{example}

\begin{remark}[Separating Points and Ensemble Controllability]
	\rm Example \ref{ex:1D_multiple} illustrates the number of independent control fields (viewed as functions of the parameter $\b$) 
	required in order to separate points of the same images on different injective branches, so that ensemble controllability can be achieved. In this example, the reachable state $\eta\in\mathcal{L}_0$ of the ensemble is of the form $\eta(\b)=\sum_kc_{1k}\b^{2k}+\sum_kc_{2k}\b^{2k+1}$, where $c_{1k},c_{2k}\in\mathbb{R}$ for $k=0,1,2,\ldots$. Although 
	$\b$ and $-\b$ have the same image under $a$, i.e., $a(\b)=a(-\b)=\b^2$, the independence of the two control vector fields, $1$ and $\b$, in $C([-1,1],\mathbb{R})$ renders ensemble controllability. Looking in a more systematic way, let $p_1(\b)=\sum_kc_{1k}\b^{2k}$ and $p_2(\b)=\sum_kc_{2k}\b^{2k}$, then the reachable states can be written as $\eta(\b)=p_1(\b)+p_2(\b)\b$. For any pair $(\bar\b,-\bar\b)$ from different injective branches of $a(\b)=\b^2$, 
	e.g., $\bar\b\in[-1,0)$ and $-\bar\b\in(0,1]$, we have $\eta(\bar\b)=h_1+h_2\bar\b$ and $\eta(-\bar\b)=h_1-h_2\bar\b$, where $h_1=p_1(\bar\b)=p_1(-\bar\b)$ and $h_2=p_2(\bar\b)=p_2(-\bar\b)$. Putting these relations into a matrix form gives $\zeta(\bar\b)=D(\bar\b)h$ with
	$$\zeta(\bar\b)=\left[\begin{array}{c} \eta(\bar\b) \\ \eta(-\bar\b)\end{array}\right], \quad D(\bar\b)=\left[\begin{array}{cc} 1 & \bar\b \\1 & -\bar\b \end{array}\right],\quad h=\left[\begin{array}{c} h_1\\ h_2 \end{array}\right].$$
	Then, if $D(\bar\b)$ is full rank for each $\bar\b\in[-1,0)\cup (0,1]$, 
	this system of linear equations has a unique solution $h$ for any given $\xi(\bar\b)$ (when $\b=0$, $h=(\eta(0),\eta(0))'$ is a solution). This implies that there exist polynomials $p_1$ and $p_2$, generated respectively by the two 
	control functions $u_1$ and $u_2$, which separate 
	$\bar\b$ and $-\bar\b$. 
	Because the two injective branches $[-1,0)$ and $(0,1]$ of $a(\b)$ are distinguishable by 
	$u_1$ and $u_2$ that generate the independent fields $\{u_1,\b u_2\}$, 
	the system in \eqref{eq:1D_single} is uniformly ensemble controllable on $C([-1,1],\mathbb{R})$.
\end{remark}

Illuminated by Example \ref{ex:1D_multiple}, systems with non-injective dynamics can be uniformly ensemble controllable when multiple inputs are available. Inspired by this finding, we establish explicit controllability conditions for an ensemble of one-dimensional time-invariant linear systems. 

\begin{theorem}
    \label{thm:1D_multiple}
    Consider an ensemble of one-dimensional multi-input linear systems indexed by the parameter $\b$ varying on a compact set $K\subset\mathbb{R}$, given by
    \begin{align}
        \label{eq:1D_multiple}
        \frac{d}{dt}x(t,\b)=a(\b)x(t,\b)+\sum_{i=1}^mb_i(\b)u_i(t),
    \end{align}
	where $x(t,\cdot)\in C(K,\mathbb{R})$, $a\in C(K,\mathbb{R})$, $b_i\in C(K,\mathbb{R})$, and $u_i:[0,T]\to\mathbb{R}$ are piecewise constant for $i=1,\ldots,m$. Then, the system in \eqref{eq:1D_multiple} is uniformly ensemble controllable on $C(K,\mathbb{R})$ if and only if ${\rm span}\{b_1|_{a^{-1}(\eta)},\dots,b_m|_{a^{-1}(\eta)}\}=C(a^{-1}(\eta),\mathbb{R})$ for all $\eta\in a(K)$, where $a^{-1}(\eta)=\{\beta\in K\mid a(\beta)=\eta\}$ is the inverse image of $\eta$ and $b_i|_{a^{-1}(\eta)}:a^{-1}(\eta)\rightarrow\mathbb{R}$ is the restriction of $b_i$ to $a^{-1}(\eta)$. 
\end{theorem}

{\it Proof:} (Sufficiency): Suppose that the condition ${\rm span}\{b_1|_{a^{-1}(\eta)},\dots,b_m|_{a^{-1}(\eta)}\}=C(a^{-1}(\eta),\mathbb{R})$ for all $\eta\in a(K)$ is satisfied, 
we wish to show that the system in \eqref{eq:1D_multiple} is uniformly ensemble controllable on $C(K,\mathbb{R})$, and equivalently, the Lie algebra generated by the drift and the control vector fields, $\mathcal{L}_0=\{ax,b_i\}_{LA}=\text{span}\{a^kb_i: k=0,1,2,\dots;\, i=1,\dots,m\}$, is dense in $C(K,\mathbb{R})$ under the uniform topology. 

We first observe that the condition ${\rm span}\{b_1|_{a^{-1}(\eta)},\dots,b_m|_{a^{-1}(\eta)}\}=C(a^{-1}(\eta),\mathbb{R})$ implies that $|a^{-1}(\eta)|\leq m$ for all $\eta\in a(K)$, i.e., $a$ has at most $m$ injective  branches, where $|a^{-1}(\eta)|$ denotes the cardinality of the set $a^{-1}(\eta)$. Without loss of generality, we assume that $a$ has exactly $m$ injective branches, then there exists a partition of $K$ containing $m$ elements, say $K_1,\dots,K_m$, such that $a|_{K_i}$ is injective for every $i=1,\dots,m$. Let $\mathcal{A}$ be the subalgebra of $C(K,\mathbb{R})$ generated by $a,b_1,\dots,b_m$, then we have $\mathcal{L}_0\subset\mathcal{A}$. By Proposition \ref{prop:1D}, for each $K_i$, the system in \eqref{eq:1D_multiple} with the parameter restricted on it is uniformly ensemble controllable on $C(K_i,\mathbb{R})$, or equivalently, $\overline{\mathcal{L}_0|_{K_i}}=C(K_i,\mathbb{R})$, where $\mathcal{L}_0|_{K_i}=\text{span}\{(a^kb_j)|_{K_i}: k=0,1,2,\dots;\, j=1,\dots,m\}$. Hence, we have $\mathcal{A}|_{K_i}\subseteq\overline{\mathcal{L}_0|_{K_i}}=C(K_i,\mathbb{R})$, where $\mathcal{A}|_{K_i}=\{f|_{K_i}: f\in\mathcal{A}\}$, and it follows that $\mathcal{A}\subseteq\overline{\mathcal{L}_0}$ because $K_i$, $i=1,\dots,m$, form a partition of $K$.

Now recall that for any two distinct points $\b_1,\b_2\in K$, if $a(\b_1)\neq a(\b_2)$, then the function $a\in\mathcal{A}$ separates these two points. 
Otherwise, if $a(\b_1)=a(\b_2)=\eta$, i.e., $\b_1,\b_2\in a^{-1}(s)$, the condition ${\rm span}\{b_1|_{a^{-1}(s)},\dots,b_m|_{a^{-1}(s)}\}=C(a^{-1}(\eta),\mathbb{R})$ then implies the existence of a continuous function $\bar{\xi}$ defined on $a^{-1}(s)$ of the form $\bar{\xi}=c_1b_1|_{a^{-1}(s)}+\cdots+c_mb_m|_{a^{-1}(\eta)}$ for some $c_1,\dots,c_m\in\mathbb{R}$, such that $\bar{\xi}(\b_1)\neq\bar{\xi}(\b_2)$. It follows that the unrestricted function $\xi=c_1b_1+\cdots+c_mb_m\in\mathcal{A}$ satisfies $\xi(\beta_1)\neq\xi(\beta_2)$, which implies that $\mathcal{A}$ separates points in $K$ (see Definition \ref{def:separating_points} in Appendix \ref{appd:S-W_thm}). Then, according to the Stone-Weierstrass theorem (Theorem \ref{thm:Weierstrass} in Appendix \ref{appd:S-W_thm}), $\mathcal{A}$ is dense in $C(K,\mathbb{R})$, and thus $C(K,\mathbb{R})=\overline{\mathcal{A}}$. Together with $\mathcal{A}\subseteq\overline{\mathcal{L}_0}\subseteq C(K,\mathbb{R})$, we obtain $\overline{\mathcal{L}_0}=C(K,\mathbb{R})$, which shows uniform ensemble controllability of the system in \eqref{eq:1D_multiple} on $C(K,\mathbb{R})$.

(Necessity): Suppose that the system in \eqref{eq:1D_multiple} is uniformly ensemble controllable on $C(K,\mathbb{R})$, then so is it on 
$C(a^{-1}(\eta),\mathbb{R})$ for any $\eta\in a(K)$. Because $a|_{a^{-1}(\eta)}$ is constant, the Lie algebra generated by the drift and control vector fields, $a|_{a^{-1}(\eta)}x$ and $b_1|_{a^{-1}(\eta)}$, $\dots$, $b_m|_{a^{-1}(\eta)}$, respectively, is $\mathcal{L}_0|_{a^{-1}(\eta)}={\rm span}\{b_1|_{a^{-1}(\eta)},\dots,b_m|_{a^{-1}(\eta)}\}$. This implies $\mathcal{L}_0|_{a^{-1}(\eta)}=C(a^{-1}(\eta),\mathbb{R})$ due to ensemble controllability of the system. 
\hfill$\Box$

Notice that, for each fixed $\eta\in a(K)$, the controllability condition 
in Theorem \ref{thm:1D_multiple} implies that ${\rm dim}(C(a^{-1}(\eta),\mathbb{R}))\leq m$. Consequently, $a^{-1}(\eta)$ must be a finite set, i.e., $a(\b)$ has finitely many injective branches, and hence a space with discrete topology. It then follows that every function defined on $a^{-1}(\eta)$ is continuous and determined by $|a^{-1}(\eta)|$ points, where $|a^{-1}(\eta)|$ denotes the cardinality of $a^{-1}(\eta)$, and thus $C(a^{-1}(\eta),\mathbb{R})$ is isomorphic to $\mathbb{R}^{|a^{-1}(\eta)|}$, because one can identify each function on $a^{-1}(\eta)$ with a $|a^{-1}(\eta)|$-dimensional vector. As a result, Theorem \ref{thm:1D_multiple} provides a `finite-dimensional' method to evaluate controllability of the ensemble system in \eqref{eq:1D_multiple} defined on the `infinite-dimensional' space $C(K,\mathbb{R})$.

\begin{corollary}
\label{cor:rank}
	Consider the linear ensemble system in \eqref{eq:1D_multiple}, and let $a^{-1}(\eta)=\{\b_\eta^1,\dots,\b_\eta^{\kappa(\eta)}\}$ be the preimage of $\eta$ under the drift $a$, where $\kappa(\eta)=|a^{-1}(\eta)|$ denotes its cardinality. Then, this system is uniformly ensemble controllable on $C(K,\mathbb{R})$ if and only if the matrix $D(\eta)\in\mathbb{R}^{\kappa(\eta)\times m}$, defined by
	\begin{equation}
		\label{eq:Ds}
		D(\eta)=\left[\begin{array}{ccc} b_1(\b_\eta^1) & \cdots & b_m(\b_\eta^1) \\ \vdots & \ddots  & \vdots \\ b_1(\b_\eta^{\kappa(\eta)}) & \cdots & b_m(\b_\eta^{\kappa(\eta)}) \end{array}\right],
	\end{equation}
\end{corollary}
has full rank, i.e., rank$(D(\eta))=\kappa(\eta)\leq m$, for all $\eta\in a(K)$. In particular, we call $D(\eta)$ the \emph{Ensemble Controllability Gramian}.

{\it Proof.} Through the isomorphism between $C(a^{-1}(\eta),\mathbb{R})$ and $\mathbb{R}^{\kappa(\eta)}$ given by $f\mapsto (f(\b_\eta^1),\dots,$ $f(\b_\eta^{\kappa(\eta)}))'$, the range space of $D(\eta)$ can be identified with the vector space spanned by $b_1|_{a^{-1}(\eta)}$, $\dots$, $b_m|_{a^{-1}(\eta)}$. As a result, the necessary and sufficient controllability condition presented in Theorem \ref{thm:1D_multiple}, i.e., ${\rm span}\{b_1|_{a^{-1}(\eta)},\dots,b_m|_{a^{-1}(\eta)}\}=C(a^{-1}(\eta),\mathbb{R})=\mathbb{R}^{\kappa(\eta)}$, can be translated to the rank condition, rank$(D(\eta))=\kappa(\eta)$. \hfill$\Box$

Consequently, if the number of control inputs $m<\kappa(s)$, 
then the system in \eqref{eq:1D_multiple} can never be uniformly ensemble controllable on $C(K,\mathbb{R})$. This observation leads to a necessary condition for uniform ensemble controllability.

\begin{remark}
    \label{rmk:1D_multiple}
    \rm
    Proposition \ref{prop:1D} showed that it requires only one control function to achieve uniform ensemble controllability for the scalar linear ensemble system with an injective drift (i.e., involving one injective branch). Theorem \ref{thm:1D_multiple}, along with Corollary \ref{cor:rank}, is a generalization to this result, illustrating that the number of independent controls must be no less than the number of injective branches of the drift, so that separating, or distinguishing, points from different injective branches with the same images is possible, leading to 
    ensemble controllability of the system.
\end{remark}

The use of the Gramian matrix to determine ensemble controllability is illuminated by the following ensemble control example, involving a system with non-injective drift driven by multiple inputs.

\begin{example}
    \label{ex:1D_multiple_1}
    \rm Consider the scalar linear ensemble system driven by two control inputs, 
		$$\frac{d}{dt}x(t,\b)=\cos(\b)x(t,\b)+u_1(t)+\beta u_2(t),$$
		where $\b\in K=[-\pi,\pi]$ is the system parameter, and $u_1,u_2:[0,T]\to\mathbb{R}$ 
		are piecewise constant 
		control functions. In this case, the drift vector field is $a(\b)=\cos(\b)$, and the inverse image of $\eta\in S$, where $S=a(K)=[-1,1]$ is the range of $a$, is given by
		$$
		a^{-1}(\eta)=\left\{\begin{array}{ll} \{0\}, & \text{if}\ \eta=1 \\
		\{\cos^{-1}(\eta),-\cos^{-1}(\eta)\}, & \text{if}\ \eta\in[-1,1),\end{array}\right.
		$$
		and hence,
	        $$
		\kappa(\eta)=\left\{\begin{array}{ll}1, & \text{if}\ \eta=1 \\
		2, & \text{if}\ \eta\in[-1,1).\end{array}\right.
		$$
		In this case, $\sup_{\eta\in S}\kappa(\eta)=2$ is equal to the number of independent controls, and the ensemble controllability Gramian, as defined in \eqref{eq:Ds},
		$$
		D(\eta)=\left\{\begin{array}{ll}\left[\begin{array}{cc} 1 & 0 \end{array}\right], & \text{if}\ \eta=1, \vspace{3pt}\\
		\left[\begin{array}{cc} 1 & \cos^{-1}(\eta) \\ 1 & -\cos^{-1}(\eta) \end{array}\right], & \text{if}\ \eta\in[-1,1), \end{array}\right.
		$$
		is of rank $\kappa(\eta)$ for all $\eta\in S$, Thus, this system is uniformly ensemble controllable on $C([-\pi,\pi],\mathbb{R})$. 
\end{example}

The controllability analysis for the one-dimensional linear ensemble system based on the notion of separating points establishes a machinery for characterizing controllability of multidimensional linear ensemble systems. In our previous work \cite{Li_TAC16}, we studied a class of finite-dimensional linear ensemble systems with linear parameter variation and showed that controllability of such ensembles is characterized by the interplay between the spectrum structure of the drift term (i.e., with or without shared components among different eigenvalues functions) and the number of independent controls. We will leverage this previous finding and adopt the concept and techniques of separating points presented above to analyze controllability of the general finite-dimensional time-invariant linear ensemble system. In particular, we will establish implementable algebraic conditions, under which the shared spectra inherited in the system dynamics can be separated in order to warrant uniform ensemble controllability.

\section{Uniform ensemble controllability of multi-dimensional linear ensemble systems}
\label{sec:n-d}
Consider the $n$-dimensional linear ensemble system,
\begin{equation}
	\label{eq:n-d}
	\frac{d}{dt}X(t,\b)=A(\b)X(t,\b)+B(\b)U(t),
\end{equation}
indexed by the parameter $\b$ taking values on a compact set $K\subset\mathbb{R}$, where $X(t,\cdot)\in C(K,\mathbb{R}^n)$ is the state, 
the control function $U:[0,T]\to\mathbb{R}^m$ is piecewise constant with $T\in (0,\infty)$, and $A\in C(K,\mathbb{R}^{n\times n})$
and $B\in C(K,\mathbb{R}^{n\times m})$ are the system and the control matrix, respectively. To avoid a trivial uncontrollable case, we, without loss of generality, assume that the control matrix $B(\b)$ has nonzero rows for any $\b\in K$ so that each state of the system in \eqref{eq:n-d} receives at least one control input. 
Our approach to analyzing controllability here is to map this ensemble of $n$-dimensional systems to an ensemble of one-dimensional systems. 
In this way, the controllability analysis 
of the system as in \eqref{eq:n-d} will directly follow the procedure developed in Section \ref{sec:1D}.

\subsection{Ensembles with diagonalizable drift}
\label{sec:diag}
We first consider the case in which $A(\b)$ is diagonalizable and has real eigenvalues for all $\b\in K$. Then, there exists a matrix $P\in C(K,{\rm GL}(n,\mathbb{R}))$ such that $A(\b)=P(\b)D(\b)P^{-1}(\b)$, where ${\rm GL}(n,\mathbb{R})$ is the group of $n\times n$ invertible real matrices and $D(\b)=\text{diag}(\lambda_1(\b),\ldots,\lambda_n(\b))$ is a diagonal matrix with $\lambda_i\in C(K,\mathbb{R})$ for each $i=1,\dots,n$. By the coordinate transformation, $Y=P^{-1}X$, the system in \eqref{eq:n-d} is transformed to 
\begin{align}
\label{eq:ensemble_linear_diagonalized}
\frac{d}{dt}Y(t,\beta)=D(\b)Y(t,\beta)+\widetilde{B}(\b)U(t), 
\end{align}
where $\widetilde{B}=P^{-1}B\in C(K,\mathbb{R}^{n\times m})$ has nonzero rows for all $\b\in K$. Because controllability of a linear system is preserved under a bijective 
transformation, it suffices to study the system in \eqref{eq:ensemble_linear_diagonalized}. In the following, we will consider two cases in terms of the injectivity of the eigenvalue functions, $\lambda_i(\b)$, of $A(\b)$.

\subsubsection{System matrix with injective eigenvalue functions} 
\label{sec:injective}
Suppose that all $\l_i\in C(K,\mathbb{R}),i=1,\dots,n$, are injective, we will derive controllability conditions for the system in \eqref{eq:n-d} by induction initialized by a two-dimensional system, given by
\begin{align}
	\label{eq:linear_ensemble_2d}
	\frac{d}{dt}\left[\begin{array}{c} y_1(t,\b) \\ y_2(t,\b) \end{array}\right]=\left[\begin{array}{cc} \l_1(\b) & 0 \\ 0 & \l_2(\b) \end{array}\right]\left[\begin{array}{c} y_1(t,\b) \\ y_2(t,\b) \end{array}\right]+\left[\begin{array}{c} \tilde{b}_1(\b) \\ \tilde{b}_2(\b) \end{array}\right]U(t),
\end{align}
where $\tilde{b}_1,\tilde{b}_2\in C(K,\mathbb{R}^{m})$ and $U:[0,T]\to\mathbb{R}^m$. Let $K_i=\l_i(K)=\{\l_i(\b):\b\in K\}$ denote the image of the function $\l_i(\b)$ for $i=1,2$, then the spectrum of $A(\b)$ is $K_1\cup K_2$. Next, we define a \emph{parameterization map} $\phi:C(K_1\cup K_2,\mathbb{R})\rightarrow C(K,\mathbb{R}^2)$ by 
$$f(\eta)\mapsto\left[\begin{array}{c} f|_{K_1}(\l_1(\b)) \\ f|_{K_2}(\l_2(\b)) \end{array}\right],$$
where $f|_{K_i}$ denotes the restriction of $f$ to $K_i$ for $i=1,2$. Then, $\phi$ is a well-defined injective function, and its image is given by 
\begin{align}
\label{eq:im_phi}
{\rm Im}(\phi)=\left\{\left[\begin{array}{c} f_1 \\ f_2 \end{array}\right]\in C(K,\mathbb{R}^2):f_1|_{\l_1^{-1}(K_1\cap K_2)}=f_2|_{\l_2^{-1}(K_1\cap K_2)}\right\}.
\end{align}

(Case I): If $K_1\cap K_2=\varnothing$, then both $f_1|_{\l_1^{-1}(K_1\cap K_2)}$ and $f_2|_{\l_2^{-1}(K_1\cap K_2)}$ are the empty function, and thus Im$(\phi)=C(K,\mathbb{R}^2)$, i.e., $\phi$ is also surjective. This follows that $\phi$ has a bijective inverse $\phi^{-1}:C(K,\mathbb{R}^2)\rightarrow C(K_1\cup K_2,\mathbb{R})$ defined by 
\begin{equation}
	\label{eq:phi_inverse}
	\left[\begin{array}{c} f_1(\b) \\ f_2(\b) \end{array}\right]\mapsto f(\eta),
\end{equation}
where 
\begin{align}
	\label{eq:f_eta}
	f(\eta)=
	\begin{cases} 
		f_1(\l_1^{-1}(\eta)),\quad\eta\in K_1 \\
		f_2(\l_2^{-1}(\eta)),\quad\eta\in K_2.
	\end{cases}
\end{align}
Applying the function 
$\phi^{-1}$ to the ensemble system in \eqref{eq:linear_ensemble_2d} and defining $z(t,\eta)=\phi^{-1}\left[\begin{array}{c} y_1(t,\b) \\ y_2(t,\b) \end{array}\right]\in\mathbb{R}$, i.e., reparameterizing the system parameterized by $\beta$ in \eqref{eq:linear_ensemble_2d} by $\eta\in K_1\cup K_2$, yields an ensemble system on $C(K_1\cup K_2,\mathbb{R})$, given by
\begin{align}
	\label{eq:linear_ensemble_nooverlapping}
	\frac{d}{dt}z(t,\eta)=\eta z(t,\eta)+b(\eta)U(t),
\end{align}
where $z(t,\cdot)\in C(K_1\cup K_2,\mathbb{R})$ and $b(\eta)=\tilde{b}_i(\l_i^{-1}(\eta))\in\mathbb{R}^{1\times m}$ for $\eta\in K_i$ and $i=1,2$.
Because $b(\eta)\neq0$ for all $\eta\in K_1\cup K_2$ by the assumption that the control matrix has no zero rows for any $\b\in K$, the scalar ensemble system in \eqref{eq:linear_ensemble_nooverlapping} is uniformly ensemble controllable on $C(K_1\cup K_2,\mathbb{R})$ according to Proposition \ref{prop:1D}, and so is the system in \eqref{eq:linear_ensemble_2d}.

(Case II): If $K_1\cap K_2\neq\varnothing$, otherwise, then ${\rm Im}(\phi)\subset C(K,\mathbb{R}^2)$, i.e., $\phi$ is not surjective anymore. However, the injectivity of $\phi$ implies the existence of a left inverse of $\phi$, say $\phi_L^{-1}:C(K,\mathbb{R}^2)\rightarrow C(K_1\cup K_2,\mathbb{R})$, and its restriction on ${\rm Im}(\phi)$, defined by $\phi_L^{-1}|_{{\rm Im}(\phi)}:{\rm Im}(\phi)\rightarrow C(K_1\cup K_2,\mathbb{R})$, is bijective. Then, applying $\phi_{L}^{-1}$ to the ensemble system in \eqref{eq:linear_ensemble_2d}, one can show that this system is ensemble controllable on ${\rm Im}(\phi)\subset C(K,\mathbb{R}^2)$. Recall \eqref{eq:im_phi} that ${\rm Im}(\phi)$ contains $\mathbb{R}^2$-valued functions restricted to $K_1\cap K_2$ with $f_1|_{\l_1^{-1}(K_1\cap K_2)}=f_2|_{\l_2^{-1}(K_1\cap K_2)}$. Therefore, 
ensemble controllability of the system in \eqref{eq:linear_ensemble_2d} on the entire state space $C(K,\mathbb{R}^2)$ will be determined by whether the points on the shared spectrum, $K_1\cap K_2$, can be separated using the available control inputs. The analysis can be carried out by considering the symmetric case in which $K_1=K_2=K\subset\mathbb{R}^+$. In this case, the parameterization function $\psi$ is defined as $\psi:C(-K\cup K,\mathbb{R})\rightarrow C(K,\mathbb{R}^2)$ by
$$f(\eta)\mapsto\left[\begin{array}{c}f|_K(\b) \\ f|_{-K}(-\b)\end{array}\right],$$
and $\psi$ is a bijection with the inverse $\psi^{-1}:C(K,\mathbb{R}^2)\rightarrow C(-K\cup K,\mathbb{R})$ defined as in \eqref{eq:phi_inverse} and \eqref{eq:f_eta}. 
The map $\psi^{-1}$ then induces 
an ensemble system on $C(-K\cup K,\mathbb{R})$, associated with the system in \eqref{eq:linear_ensemble_2d}, of the form
\begin{align}
	\label{eq:linear_ensemble_overlapping}
	\frac{d}{dt}z(t,\eta)=|\eta|\, z(t,\eta)+b(\eta)U(t),
\end{align}
where $b(\eta)=\tilde{b}_1(\l_1^{-1}(\eta))$ for $\eta\in K$, and $b(\eta)=\tilde{b}_2(\l_2^{-1}(-\eta))$ for $\eta\in -K$. Because the drift $\varphi(\eta)=|\eta|$
has two injective branches $-K$ and $K$, by Corollary \ref{cor:rank}, the system in \eqref{eq:linear_ensemble_overlapping} is uniformly ensemble controllable on $C(-K\cup K,\mathbb{R})$ if and only if the ensemble controllability Gramian $D(\eta)\in\mathbb{R}^{2\times m}$ has rank 2 for all $\eta\in K$, where
\begin{equation}
	\label{eq:D_eta}
	D(\eta)=\left[\begin{array}{ccc} b_1(\eta) & \cdots & b_m(\eta) \\ b_1(-\eta) & \cdots & b_m(-\eta) \end{array}\right]=\left[\begin{array}{c} b(\eta)  \\ b(-\eta) \end{array}\right].
\end{equation}

\begin{figure}[t]
     \centering
     \includegraphics[width=1\columnwidth]{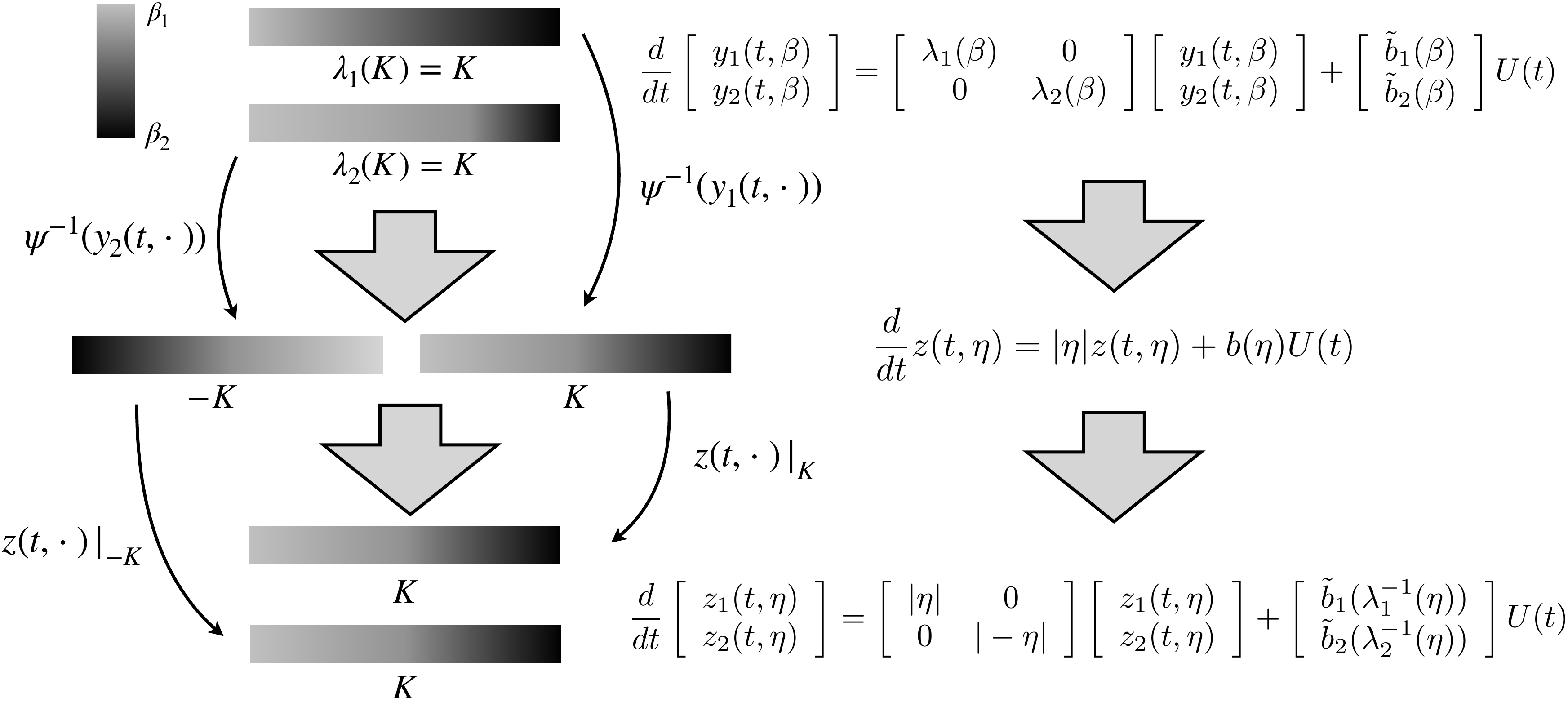}
     \caption{A schematic illustration of the reparameterization process for separating points in the overlapping spectrum of the system in \eqref{eq:linear_ensemble_2d}. This figure shows the symmetric case where $\l_1(K)=\l_2(K)=K=[\b_1,\b_2]\subset\mathbb{R}^+$, in which the parameterization map $\psi$ transforms the two-dimensional ensemble system $Y(t,\b)$ to a scalar ensemble $z(t,\eta)$ with the drift $\varphi(\eta)=|\eta|$ that has two injective branches $K$ and $-K$. The dynamics of the scalar ensemble restricted to $K$ and $-K$ can be represented using two states over distinct domains, i.e., $z_1(t,\cdot)\doteq z(t,\cdot)|_K$ and $z_2(t,\cdot)\doteq z(t,\cdot)|_{-K}$, which in this case then constitute the two-dimensional system reparameterized by $\eta$ (see \eqref{eq:linear_ensemble_overlapping_reduced_2} as well). 
     The different gray levels in the bars denote different values in the continuous spectrum, $\l_1(\b)$, $\l_2(\b)$, $\varphi(\eta)$, $\varphi|_K(\eta)$, or $\varphi|_{-K}(\eta)$, of the system matrix of the corresponding ensemble on the right hand side.}
     \label{fig:overlapping}
\end{figure}

From the perspective of separating points, rank$(D(\eta))=2$ indicates that the two injective branches of the drift $|\eta|$ in system \eqref{eq:linear_ensemble_overlapping} are separable by the available control inputs. In other words, in this case, \emph{this one-dimensional system behaves like a two-dimensional system with each state representing the dynamics of the system in \eqref{eq:linear_ensemble_overlapping} restricted to one of the two injective branches of $|\eta|$,} which can be expressed as
\begin{align}
	\label{eq:linear_ensemble_overlapping_reduced_1}
	\frac{d}{dt}\left[\begin{array}{c} z_1(t,\eta_1) \\ z_2(t,\eta_2) \end{array}\right]
	&= \left[\begin{array}{cc} \eta_1 & 0 \\ 0 & \eta_2 \end{array}\right]\left[\begin{array}{c} z_1(t,\eta_1) \\ z_2(t,\eta_2) \end{array}\right]+\left[\begin{array}{c} \tilde{b}_1(\l_1^{-1}(\eta_1)) \\  \tilde{b}_2(\l_2^{-1}(\eta_2))\end{array}\right]U(t)
\end{align}
for $(\eta_1,\eta_2)=(|\eta|,|-\eta|)\in K\times K$. Therefore, in this case, we have $\eta_1=\eta_2=\eta\in K$ so that
\begin{align}
	\frac{d}{dt}\left[\begin{array}{c} z_1(t,\eta) \\ z_2(t,\eta) \end{array}\right] &= \left[\begin{array}{cc} \eta & 0 \\ 0 & \eta \end{array}\right]\left[\begin{array}{c} z_1(t,\eta) \\ z_2(t,\eta) \end{array}\right]+\left[\begin{array}{c} b(\eta) \\  b(-\eta)\end{array}\right]U(t), \nonumber\\
	\label{eq:linear_ensemble_overlapping_reduced_2}
	&= \eta I_2\left[\begin{array}{c} z_1(t,\eta) \\ z_2(t,\eta) \end{array}\right]+D(\eta)U(t),
\end{align}
where $I_2$ is the 2-by-2 identity matrix. A schematic illustration of transforming the system in \eqref{eq:linear_ensemble_2d} to \eqref{eq:linear_ensemble_overlapping_reduced_2} by the parameterization $\eta_i=\lambda_i(\b)$ such that $z_i(t,\eta)=y_i(t,\l_i^{-1}(\eta))$, $i=1,2$, is presented in Figure \ref{fig:overlapping}. Because $\l_1$ and $\l_2$ are injective, their inverse functions $\l_1^{-1}$ and $\l_2^{-1}$ exist so that such parameterizations are always well-defined regardless of whether $K_1$ and $K_2$ are disjoint or overlapped.

Observe that the ``classical'' controllability Gramian matrix associated with each individual system $\eta$ 
in \eqref{eq:linear_ensemble_overlapping_reduced_1} is given by
$$W=\left[\begin{array}{cc} \tilde{b}_1(\l_1^{-1}(\eta)) & \eta\tilde{b}_1(\l_1^{-1}(\eta)) \\  \tilde{b}_2(\l_2^{-1}(\eta)) & \eta\tilde{b}_2(\l_2^{-1}(\eta))\end{array}\right]=\big[D(\eta) \ \big|\ \eta D(\eta)\big]\in\mathbb{R}^{2\times 2m},$$
which has the same rank as the ``ensemble'' controllability Gramian, i.e., ${\rm rank}(W)={\rm rank}(D(\eta))=2$ for all $\eta\in K$. This equivalence gives rise to a characterization of the ensemble system by examining controllability of each individual system in the ensemble. For example, the ensemble system in \eqref{eq:linear_ensemble_2d} is ensemble controllable on $C(K_1\cup K_2,\mathbb{R})$ if and only if the two-dimensional system in \eqref{eq:linear_ensemble_overlapping_reduced_1} is controllable on $\mathbb{R}^2$ for all $(\eta_1,\eta_2)\in K_1\times K_2$, which will be shown below in Lemma \ref{lem:linear_2d}. Note that in general controllability of each individual system in the ensemble does not imply ensemble controllability.

\begin{lemma}
	\label{lem:linear_2d}
	The ensemble system in \eqref{eq:linear_ensemble_2d} is uniformly ensemble controllable on $C(K,\mathbb{R}^2)$ if and only if each individual system, characterized by $(\eta_1,\eta_2)\in K_1\times K_2$, in the induced ensemble system parameterized by $\eta_1=\l_1(\b)$ and $\eta_2=\l_2(\b)$ in \eqref{eq:linear_ensemble_overlapping_reduced_1} is controllable on $\mathbb{R}^2$, where $K_1=\l_1(K)$ and $K_2=\l_2(K)$.
\end{lemma}

{\it Proof.}
The proof follows from the existence of a bijective parameterization $\psi: C(K',\mathbb{R})\rightarrow C(K,\mathbb{R}^2)$ (see Appendix \ref{app:2-d}) such that the dynamics of  $z(t,\eta)\doteq\, \psi^{-1}\left[\begin{array}{c} y_1(t,\b) \\ y_2(t,\b) \end{array}\right]$ obey the scalar linear ensemble system
\begin{align}
	\label{eq:1-d-param}
	\frac{d}{dt}z(t,\eta)=\varphi(\eta)z(t,\eta)+b(\eta)U(t),
\end{align}
where $K'=K_1'\cup^* K_2'$ is the disjoint union of the two sets $K_1'$ and $K_2'$ that are homeomorphic to $K_1$ and $K_2$, respectively; $\varphi:K'\rightarrow K_1\cup K_2$ is continuous with $\varphi^{-1}(K_i)=K_i'$ for $i=1,2$, satisfying 
$$|\varphi^{-1}(s)|=\left\{\begin{array}{ll} 2, & \text{for}\ s\in K_1\cap K_2 \ \text{and}\ K_1\cap K_2\neq\varnothing,\\
1, & \text{for}\ s\in K_1\backslash K_2\ \text{or}\ s\in K_2\backslash K_1; \end{array}\right.
$$
and the control vector field 
$$b(\eta)=\left\{\begin{array}{ll}\tilde{b}_1(\l_1^{-1}\circ\varphi(\eta)), & \text{for}\ \eta\in K_1', \\
\tilde{b}_2(\l_2^{-1}\circ\varphi(\eta)), & \text{for}\ \eta\in K_2'. \end{array}\right.$$
In this case, the original ensemble of two-dimensional systems in \eqref{eq:linear_ensemble_2d} is uniformly ensemble controllable on $C(K,\mathbb{R}^2)$ if and only if the reparameterized ensemble of one-dimensional systems in \eqref{eq:1-d-param} has the same property on $C(K',\mathbb{R})$.

By Corollary \ref{cor:rank}, ensemble controllability of the system as in \eqref{eq:1-d-param} can be evaluated by the rank of the ensemble controllability Gramian,
\begin{align*}
D(s)=
\begin{cases}
\ \tilde{b}_1(\l_1^{-1}(s)) ,\qquad\quad\text{if } s\in K_1\backslash K_2,\\
\ \tilde{b}_2(\l_2^{-1}(s)),\qquad\quad\text{if } s\in K_2\backslash K_1,\\
\ \left[\begin{array}{c} \tilde{b}_1(\l_1^{-1}(s)) \\  \tilde{b}_2(\l_2^{-1}(s))\end{array}\right],\quad\text{if } s\in K_1\cap K_2.
\end{cases}
\end{align*}
Therefore, it suffices to prove that $D(s)$ is full rank for any $s\in K_1\cup K_2$ if and only if the system in \eqref{eq:linear_ensemble_overlapping_reduced_1} is controllable on $\mathbb{R}^2$ for any $(\eta_1,\eta_2)\in K_1\times K_2$.

Because $\tilde{B}$ has nonzero rows, ${\rm rank}(\tilde{b}_1(\lambda_1^{-1}((s)))=1$ and ${\rm rank}(\tilde{b}_2(\lambda_2^{-1}((s)))=1$ hold for any $s\in K_1$ and $s\in K_2$, respectively, which implies the full rank of $D(s)$ for $s\in K_1\backslash K_2$ and $s\in K_2\backslash K_1$. Correspondingly, the system in \eqref{eq:linear_ensemble_overlapping_reduced_1} is controllable on $\mathbb{R}^2$ for $\eta_1\neq\eta_2$ because in this case the controllability Gramian, for each $\eta=(\eta_1,\eta_2)'$,
$$W(\eta)=\left[\begin{array}{cc} \tilde{b}_1(\l_1^{-1}(\eta_1)) & \eta_1\tilde{b}_1(\l_1^{-1}(\eta_1)) \\  \tilde{b}_2(\l_2^{-1}(\eta_2)) & \eta_2\tilde{b}_2(\l_2^{-1}(\eta_2)) \end{array}\right],$$
is of rank 2. On the other hand, for the case in which $\eta_1=\eta_2=s\in K_1\cap K_2$, we have $W(s)=\big[D(s) \ \big|\ sD(s)\big]$ with ${\rm rank}(W(s))={\rm rank}(D(s))$. Therefore, the system in \eqref{eq:linear_ensemble_overlapping_reduced_1} is controllable on $\mathbb{R}^2$ if and only if $D(s)$ is full rank, i.e., ${\rm rank}(W(s))={\rm rank}(D(s))=2$. \hfill$\Box$

Following the same argument, the result in Lemma \ref{lem:linear_2d} can be directly extended and applied to the ensemble of $n$-dimensional linear systems defined on $C(K,\mathbb{R}^n)$.

\begin{theorem}
	\label{thm:linear_diagonalizable_injective}
	Consider the time-invariant linear ensemble system,
	\begin{align}
		\label{eq:X(t,b)}
		\frac{d}{dt}X(t,\b)=A(\b)X(t,\b)+B(\b)U(t),
	\end{align}	
	indexed by the parameter $\b$ taking values on a compact set $K\subset\mathbb{R}$, where $X(t,\cdot)\in C(K,\mathbb{R}^n)$, $U:[0,T]\rightarrow\mathbb{R}^m$ is piecewise constant, $A\in C(K,\mathbb{R}^{n\times n})$ is diagonalizable with real eigenvalues, and $B\in C(K,\mathbb{R}^{n\times m})$. Let 
	\begin{align}
		\label{eq:Y(t,b)}
		\frac{d}{dt}Y(t,\beta)=\Lambda(\b)Y(t,\beta)+\widetilde{B}(\b)U(t)
	\end{align}	
	be the corresponding diagonalized system, transformed by the eigenvalue decomposition, where 
	$$Y(t,\b)=\left[\begin{array}{c} y_1(t,\b) \\ \vdots \\ y_n(t,\b) \end{array}\right], \quad \Lambda(\b)=\left[\begin{array}{ccc} \l_1(\b) &  &   \\ & \ddots & \\ &   & \l_n(\b) \end{array}\right], \quad \widetilde{B}(\b)
	=\left[\begin{array}{c} \tilde{b}_1(\b) \\ \vdots \\ \tilde{b}_n(\b) \end{array}\right],$$
	in which $\l_i\in C(K,\mathbb{R})$, $i=1,\ldots,n$, are injective and $\tilde{b}_i\in C(K,\mathbb{R}^m)$ is the $i^{\rm th}$ row of $\widetilde{B}$. This system is uniformly ensemble controllable on $C(K,\mathbb{R}^n)$ if and only if the induced system of $Y(t,\b)$ parameterized by $\eta_1=\l_1(\b)$, $\dots$, $\eta_n=\l_n(\b)$, given by
	\begin{align}
		\label{eq:linear_ensemble_reparameterized}
		\frac{d}{dt}\left[\begin{array}{c} z_1(t,\eta_1) \\ \vdots \\ z_n(t,\eta_n) \end{array}\right]=\left[\begin{array}{ccc} \eta_1 &  &   \\ & \ddots & \\ &   & \eta_n \end{array}\right]\left[\begin{array}{c} z_1(t,\eta_1) \\ \vdots \\ z_n(t,\eta_n) \end{array}\right]+\left[\begin{array}{c} \tilde{b}_1(\l_1^{-1}(\eta_1)) \\ \vdots \\ \tilde{b}_n(\l_n^{-1}(\eta_n)) \end{array}\right]U(t),
	\end{align}
	is controllable on $\mathbb{R}^n$ for each $n$-tuple $(\eta_1,\dots,\eta_n)\in K_1\times\cdots\times K_n$, where $K_i=\l_i(K)$ for $i=1,\dots,n$. 
\end{theorem}

{\it Proof.}
We will prove the theorem by induction on the system dimension $n$. The base case of $n=2$ was shown in Lemma \ref{lem:linear_2d}. Now, suppose that the theorem holds for $n=k$, $k\geq 3$. For $n=k+1$, without loss of generality, we assume that the theorem holds for the subsystem of the ensemble 
$Y(t,\b)$ containing the first $k$ states. Then, similar to the case of $n=2$ presented in Lemma \ref{lem:linear_2d}, there exists a bijective function $\psi: C(K_1'\cup\cdots\cup K_k',\mathbb{R})\rightarrow C(K_1\cup\cdots\cup K_k,\mathbb{R}^k)$, which parameterizes the first $k$ states $y_1(t,\b),\ldots,y_k(t,\b)$ of $Y(t,\b)$ to an ensemble system on $C(K_1'\cup\cdots\cup K_k',\mathbb{R})$ of the form
\begin{align*}
\frac{d}{dt}z(t,\eta)=\varphi(\eta)z(t,\eta)+b(\eta)U(t),
\end{align*}
where $K_i'$ are disjoint compact subsets of $\mathbb{R}$ that are homeomorphic to $K_i$ for each $i=1,\dots,k$, $\varphi:K_1'\cup\cdots\cup K_k'\rightarrow K_1\cup\cdots\cup K_k$ is continuous with $\varphi^{-1}(K_i)=K_i'$, and $b(\eta)=\tilde{b}_i(\l_i^{-1}\circ\varphi(\eta))$ for $\eta\in K_i'$. Then, applying $\varphi$ to the system $Y(t,\b)$ results in an ensemble of two-dimensional linear systems parameterized by $(\eta,\eta_{k+1})\in(K_1'\cup\cdots\cup K_{k}')\times K_{k+1}$,
\begin{equation*}
	\frac{d}{dt}\left[\begin{array}{c} z(t,\eta) \\ z_{k+1}(t,\eta_{k+1}) \end{array}\right]= \left[\begin{array}{cc} \varphi(\eta) & 0 \\ 0 & \eta_{k+1} \end{array}\right]\left[\begin{array}{c} z(t,\eta) \\ z_{k+1}(t,\eta_{k+1}) \end{array}\right]\\
+\left[\begin{array}{c} b(\eta) \\  \tilde{b}_{k+1}({\l_{k+1}^{-1}(\eta_{k+1})})\end{array}\right]U,
\end{equation*}
and then the theorem follows from the induction hypothesis. 
\hfill$\Box$

\begin{remark}
\rm
In Theorem \ref{thm:linear_diagonalizable_injective}, we established a nontrivial connection of inferring uniform ensemble controllability of a linear ensemble system by controllability of each individual system in a reparameteriezed form. However, this theorem is by no means to indicate that ensemble controllability of an ensemble system on $C(K,\mathbb{R}^n)$ is implied by controllability of each individual system in the ensemble on $\mathbb{R}^n$. It is though necessary that if the ensemble is uniformly ensemble controllable on $C(K,\mathbb{R}^n)$, then each individual system must be controllable on $\mathbb{R}^n$, but the converse is not true. Note that if we treat the system in \eqref{eq:linear_ensemble_reparameterized} as an ensemble system defined on $C(K_1\times\cdots K_n,\mathbb{R}^n)$, then it is an equivalent ensemble system to $Y(t,\b)$ and to $X(t,\b)$, so that their controllability property remains the same.
\end{remark}

\subsubsection{System matrix with non-injective eigenvalue functions} 
\label{sec:non-injective}
As opposed to Section \ref{sec:injective}, in this section we analyze the situation where the system matrix of the ensemble has non-injective eigenvalue functions. 
In this case, the ability to separate ``points of the same image from different injective branches'' and ``points in each overlapping spectrum'' is the key to achieve ensemble controllability, and hence we leverage the respective results developed in Corollary \ref{cor:rank} (related to separating points in each injective branch) and Theorem \ref{thm:linear_diagonalizable_injective} (related to separating points in each overlapping spectrum) to establish ensemble controllability conditions.

\begin{theorem}
	\label{thm:linear_ensemble}
	Consider the time-invariant linear ensemble system $X(t,\b)$ as in \eqref{eq:X(t,b)} and its diagonalized ensemble system $Y(t,\b)$ as in \eqref{eq:Y(t,b)}, in which $\l_i\in C(K,\mathbb{R})$, $i=1,\ldots,n$, and $\tilde{b}_i=(\tilde{b}_{i1},\ldots,\tilde{b}_{im})\in C(K,\mathbb{R}^m)$ is the $i^{\rm th}$ row of $\widetilde{B}$. 
	This system is uniformly ensemble controllable on $C(K,\mathbb{R}^n)$ if and only if the induced system of $Y(t,\b)$ parameterized by $\eta_1=\l_1(\b),\ldots,\eta_n=\l_n(\b)$, 
	given by
	\begin{align*}
		\frac{d}{dt}\left[\begin{array}{c} Z_1(t,\eta_1) \\ \vdots \\ Z_n(t,\eta_n) \end{array}\right]=\left[\begin{array}{ccc} \eta_1I_{\kappa_1(\eta_1)} &  &   \\ & \ddots & \\ &   & \eta_nI_{\kappa_n(\eta_n)} \end{array}\right]\left[\begin{array}{c} Z_1(t,\eta_1) \\ \vdots \\ Z_n(t,\eta_n) \end{array}\right]+\left[\begin{array}{c} D_1(\eta_1) \\ \vdots \\ D_n(\eta_n) \end{array}\right]U(t), 
	\end{align*}
	is controllable on $\mathbb{R}^N$ for each $n$-tuple $(\eta_1,\dots,\eta_n)\in K_1\times\cdots\times K_n$ with $K_i=\l_i(K)$, $i=1,\ldots,n$, where $\kappa_i(\eta_i)=|\l_i^{-1}(\eta_i)|$ is the cardinality of the preimage of 
	$\eta_i$ under $\lambda_i$, $N=\sum_{i=1}^n\kappa_i(\eta_i)$, $I_{\kappa_i(\eta_i)}$ is the $\kappa_i(\eta_i)\times\kappa_i(\eta_i)$ identity matrix, and $D_i(\eta_i)\in\mathbb{R}^{\kappa_i(\eta_i)\times m}$ is the ensemble controllability Gramian associated with the $i^{\rm th}$ state variable of the ensemble $Y(t,\b)$, 
	i.e., $\frac{d}{dt}y_i(t,\b)=\l_i(\b)y_i(t,\b)+\tilde{b}_i(\b)U(t)$. 
\end{theorem}

{\it Proof.}
	By Corollary \ref{cor:rank}, each sub-ensemble $y_i(t,\b)$ is ensemble controllable on $C(K,\mathbb{R})$ if and only if the ensemble controllability Gramian 
	$$D_i(\eta_i)=\left[\begin{array}{c} \tilde{b}_{i}(\beta_{\eta_i}^1) \\ \hline \vdots \\ \hline  \tilde{b}_{i}(\beta_{\eta_i}^{\kappa(\eta_i)}) \end{array}\right]\in\mathbb{R}^{\kappa_i(\eta_i)\times m}$$ 
	is of rank $\kappa_i(\eta_i)$ for every $\eta_i\in K_i$. This implies that for each $\eta_i\in K_i\subset\mathbb{R}$, the 
	system 
	\begin{align}
		\label{eq:non_injective}
		\frac{d}{dt}Z_i(\eta_{i})=\eta_iI_{\kappa_i(\eta_i)}Z_i(t,\eta_i)+D_i(\eta_i)U(t)
	\end{align}
	is controllable on $\mathbb{R}^{\kappa_i(\eta_i)}$, because the controllability Gramian of the system in \eqref{eq:non_injective} $W_i(\eta_i)=[D(\eta_i)\mid\eta_iD_i(\eta_i)\mid\cdots\mid\eta_i^{\kappa_i(\eta_i)-1}D_i(\eta_i)]$ has the same rank as $D_i(\eta_i)$. 
	Since each state variable of the system in \eqref{eq:non_injective} represents the dynamics of $y_i(t,\b)$ 
	restricted to one injective branch of $\l_i$, 
	 the rest of the proof directly follows the same case discussed in Theorem \ref{thm:linear_diagonalizable_injective}.
	\hfill$\Box$

\begin{example}\rm
	Consider the linear ensemble system with linear parameter variation in the system matrix, given by
	\begin{align}
		\label{eq:bA}
		\frac{d}{dt}X(t,\b)=\b AX(t,\beta)+BU(t),
	\end{align}
	where $X(t,\cdot)\in C(K,\mathbb{R}^n)$ is the state, $\b\in K$ is the parameter, $K$ is a compact subset of $\mathbb{R}$, $A\in\mathbb{R}^{n\times n}$ is diagonalizable with real eigenvalues, $B\in\mathbb{R}^{n\times m}$ has nonzero rows, and $U:[0,T]\rightarrow\mathbb{R}^m$ is piecewise constant. In our previous work \cite{Li_TAC16}, ensemble controllability conditions of this system were derived 
	using the technique of polynomial approximation, and here we will retreat these conditions through the lens of separating points. 

	Without loss of generality, it suffices to assume that $A$ is a diagonal matrix, i.e., $A={\rm diag}(\rho_1,\dots,\rho_n)$, with $\rho_i\in\mathbb{R}$ for all $i=1,\dots,n$. We first argue that $\rho_i\neq0$ 
	is necessary to guarantee ensemble controllability for this system,	which implies 
	${\rm rank}(A)=n$. By contradiction, if $\rho_i=0$ for some $i$, then the $i^{\rm th}$ eigenvalue 
	of the system matrix $\b A$ is a constant function, that is, $\l_i(\b)=\beta\rho_i=0$. Hence, $\l_i^{-1}(0)=K$ contains infinitely many points so that $\l_i$ has infinitely many injective branches, which cannot be separated using a finite number of control inputs (see Remark \ref{rmk:1D_multiple}). As a result, the $i^{\rm th}$ state variable 
	is not ensemble controllable on $C(K,\mathbb{R})$, namely, 
	the system in \eqref{eq:bA} is not ensemble controllable on $C(K,\mathbb{R}^n)$.


	Next, 
	because $\l_i(\beta)=\beta\rho_i$ is injective when $\rho_i\neq0$, by Theorem \ref{thm:linear_diagonalizable_injective}, the system in \eqref{eq:bA} is ensemble controllable on $C(K,\mathbb{R}^n)$ if and only if its reparameterized form
	\begin{align}
		\label{eq:bA_repara}
		\frac{d}{dt}Z(t,\eta)=\Lambda(\eta)Z(t,\eta)+BU(t)
	\end{align}
	is controllable on $\mathbb{R}^n$ for every $\eta=(\eta_1,\dots,\eta_n)\in K_1\times\cdots\times K_n$, where 
	$K_i=\l_i(K)=\{\eta_i:\eta_i=\beta\rho_i\text{ for }\beta\in K\}$ and $\Lambda(\eta)={\rm diag}(\eta_1,\dots,\eta_n)$. To analyze controllability of the system in \eqref{eq:bA_repara}, we consider the following two cases:

	(i) If $0\in K$, then we have $0\in K_i$ for every $i=1,\dots,n$. It follows that $n$ independent controls are required, i.e., ${\rm rank}(B)=n$, to separate the shared eigenvalue $0$ for the system in \eqref{eq:bA_repara}, and equivalently in \eqref{eq:bA}, to be ensemble controllable. Alternatively, this is because that the subsystem with $\eta=(0,\dots,0)$, i.e., $\frac{d}{dt}Z(t,0)=BU(t)$, is controllable when ${\rm rank}(B)=n$.

	(ii) Suppose that $0\not\in K$ and that $\xi$ is in a shared spectrum among $K_i$, $i=1,\ldots,n$. Without loss of generality, here we assume that $\xi\in\cap_{j=1}^r K_j$. Then, \eqref{eq:bA_repara} becomes
	$$\frac{d}{dt}\left[\begin{array}{c} z_1(t,\eta) \\ \vdots \\ z_r(t,\eta) \\  z_{r+1}(t,\eta_{r+1}) \\ \vdots \\ z_n(t,\eta_n) \end{array}\right]=\left[\begin{array}{cccccc} \xi & & & & &  \\ & \ddots & & & & \\ &  & \xi & & &  \\  &  & & \eta_{r+1} & & \\ & & & & \ddots & \\ & & & & & \eta_n \end{array}\right]\left[\begin{array}{c} z_1(t,\eta) \\ \vdots \\ z_r(t,\eta) \\  z_{r+1}(t,\eta_{r+1}) \\ \vdots \\ z_n(t,\eta_n) \end{array}\right]+\left[\begin{array}{c} {b}_1 \\ \vdots \\  {b}_r \\ {b}_{r+1} \\ \vdots \\ {b}_n \end{array}\right]U(t),$$
	where $b_i$, $i=1,\ldots,n$, is the $i^{th}$ row of $B$. Controllability of the above system requires separating the shared eigenvalue $\xi$ among the first $r$ states $z_1$, $\cdots$, $z_r$, and, equivalently, requires the full rank of the ensemble controllability Gramian associated with the subsystem containing these $r$ state variables, i.e., $\text{rank}(W_r)=r$, where
	$$W_r=\left[\begin{array}{c|c|c|c} b_1 & \xi b_1 & \cdots & \xi^{r-1} b_1 \\ \vdots & \vdots & \cdots & \vdots \\ b_r & \xi b_r & \cdots & \xi^{r-1} b_r \end{array}\right]=\Big[B_r\mid\xi B_r\mid\cdots\mid \xi^{r-1}B_r\Big]$$ 
	and $B_r=(b_1',\ldots,b_r')'$ is an $r\times m$-submatrix of $B$. Because ${\rm rank}(B_r)={\rm rank}(W_r)=r$, this shows that $r$ independent controls applied to the first $r$ states are necessary to separate the point $\xi$ and render uniform ensemble controllability.

	The ensemble controllability conditions established above in (i) and (ii) for the system in \eqref{eq:bA} by using the method of separating points coincide with the conditions presented in our previous work \cite{Li_TAC16}.
\end{example}

\subsection{Ensembles with non-diagonalizable drift}
\label{sec:nondiagonalizable}
In this section, we take a step further to analyze controllability of linear ensembles whose system matrices are not diagonalizable. We start off our investigation by considering the ensemble system with its system matrix being a single Jordan block, given by
\begin{align*}
\frac{d}{dt}X(t,\b)=J(\b)X(t,\b)+B(\b)U(t),
\end{align*}
where $\b$ is the parameter varying on a compact set $K\subset\mathbb{R}$, 
$$J(\b)=\left[\begin{array}{ccccc} \lambda(\b) & 1 &  & & \\  & \lambda(\b) & 1 & & \\  &  & \ddots & \ddots &  \\  & & & \ddots & 1\\ &  &  & & \lambda(\b) \end{array}\right],$$
and $\l\in C(K,\mathbb{R})$ is injective. To fix ideas, let's consider the two-dimensional case,
\begin{align}
	\frac{d}{dt}\left[\begin{array}{c} x_1(t,\b) \\ x_2(t,\b) \end{array}\right]&= \left[\begin{array}{cc} \lambda(\b) & 1 \\ 0 &  \lambda(\b) \end{array}\right]\left[\begin{array}{c} x_1(t,\b) \\ x_2(t,\b) \end{array}\right]+\left[\begin{array}{c} b_1(\b) \\ b_2(\b) \end{array}\right]U(t),\nonumber\\
	&=\left[\begin{array}{cc} \lambda(\b) & 0 \\ 0 &  \lambda(\b) \end{array}\right]\left[\begin{array}{c} x_1(t,\b) \\ x_2(t,\b) \end{array}\right]+\left[\begin{array}{c} b_1(\b)U(t)+x_2(t,\b) \\ b_2(\b)U(t) \end{array}\right],\label{eq:jordan_2d}
\end{align}
where $b_i(\b)$ is the $i^{\rm th}$ row of $B(\b)$, and we assume that $b_i(\b)\neq0$ for $i=1,2$ and for all $\b\in K$. Since $x_2(t,\b)$ depends on $U(t)$, 
we can write an equivalent control vector field for the system \eqref{eq:jordan_2d}, 
that is,
\begin{align*}
\hat{B}(\b) &= \frac{\partial}{\partial U(t)}\left[\begin{array}{c} b_1(\b)U(t)+x_2(t,\b) \\ b_2(\b)U(t) \end{array}\right]=\left[\begin{array}{c} b_1(\b)+\frac{\partial}{\partial U(t)}x_2(t,\b) \\ b_2(\b) \end{array}\right], \\
&= \left[\begin{array}{c} b_1(\b)+tb_2(\b) \\ b_2(\b) \end{array}\right]\in\mathbb{R}^{2\times m},
\end{align*}
where the last equality is due to the relation $\frac{d}{dt}\frac{\partial}{\partial U(t)}x_2(t,\b)=\frac{\partial}{\partial U(t)}\frac{d}{dt}x_2(t,\b)=b_2(\b)$, 
which holds piecewisely in $t$ because $U(t)$ is piecewise constant, so that $\frac{\partial x_2(t,\b)}{\partial U(t)}=tb_2(\b)$. As a result, $\text{rank}(\hat{B}(\b))=\text{rank}(B(\b))$ for all $\b\in K$, so that the controllability characterization of the system in \eqref{eq:jordan_2d} is equivalent to its diagonalizable counterpart,
\begin{align}
	\label{eq:diag_2d}
	\frac{d}{dt}\left[\begin{array}{c} x_1(t,\b) \\ x_2(t,\b) \end{array}\right]=\left[\begin{array}{cc} \lambda(\b) & 0 \\ 0 &  \lambda(\b) \end{array}\right]\left[\begin{array}{c} x_1(t,\b) \\ x_2(t,\b) \end{array}\right]+\left[\begin{array}{c} b_1(\b) \\ b_2(\b) \end{array}\right]U(t).
\end{align}
For example, if $\text{rank}(\hat{B}(\b))=\text{rank}(B(\b))=2$, by Theorem \ref{thm:linear_diagonalizable_injective}, the system in \eqref{eq:diag_2d} is uniformly ensemble controllable, and so is the system is \eqref{eq:jordan_2d}.

The above analysis can be directly generalized to the higher dimensional case, described in the following theorem.

\begin{proposition}
	\label{prop:jordan_block}
	Consider the time-invariant linear ensemble system, $\frac{d}{dt}X(t,\b)$ $=$ $A(\b)X(t,\b)+B(\b)U(t)$, 
	indexed by the parameter $\b$ taking values on a compact set $K\subset\mathbb{R}$, where $X(t,\cdot)\in C(K,\mathbb{R}^n)$, $U:[0,T]\rightarrow\mathbb{R}^m$ is piecewise constant, $A\in C(K,\mathbb{R}^{n\times n})$ is similar to a Jordan block,
	$$J(\b)=\left[\begin{array}{ccccc} \lambda(\b) & 1 &  & & \\  & \lambda(\b) & 1 & & \\  &  & \ddots & \ddots &  \\  & & & \ddots & 1\\ &  &  & & \lambda(\b) \end{array}\right],$$
	with $\l\in C(K,\mathbb{R})$ injective, and $B\in C(K,\mathbb{R}^{n\times m})$. This system is uniformly ensemble controllable on $C(K,\mathbb{R}^n)$ if and only if rank$(B(\b))=n$ for all $\b\in K$. 
\end{proposition}

{\it Proof.}
The proof follows by the induction on $n$. The case of $n=2$ was analyzed above, and we assume that the theorem holds for $n=k$ for $k>2$. In the case that $n=k+1$, consider the subsystem $X_k(t,\b)$ consisting of $k$ states of $X(t,\b)$, given by $\frac{d}{dt}X_k(t,\b)=J_k(\b)X_k(t,\b)+B_k(\b)U(t)$, where $X_k(t,\b)=(x_2(t,\b),\dots,x_{k+1}(t,\beta))'$ and $J_k(\b)\in\mathbb{R}^{k\times k}$ and $B_k(\b)\in\mathbb{R}^{k\times m}$ are the corresponding system and control matrices. Then, the induction hypothesis implies that the subsystem $X_k(t,\b)$ is uniformly ensemble controllable on $C(K,\mathbb{R}^k)$ if and only if rank$(B_k(\b))=k$ for all $\b\in K$. Equivalently, for each fixed $\b\in K$, the range space of $B_k(\b)$, denoted by ${\rm Im}(B_k(\b))$, is isomorphic to $\mathbb{R}^k$, and then the ensemble system $X(t,\b)$ can be written as
\begin{align}
\label{eq:jordan_n_2}
\frac{d}{dt}\left[\begin{array}{c} x_1(t,\b) \\ X_k(t,\b) \end{array}\right]=\left[\begin{array}{cc} \l(\b) & 1 \\   & J_k(\b) \end{array}\right]\left[\begin{array}{c} x_1(t,\b) \\ X_k(t,\b) \end{array}\right]+\left[\begin{array}{c} b_1(\b) \\ B_k(\b) \end{array}\right]U(t).
\end{align}
This system can be treated as a 
linear ensemble system defined on $C(K,\mathbb{R})\oplus C(K,\mathbb{R}^k)$ with the system matrix similar to a Jordan block, where $\oplus$ denotes the direct sum of vector spaces. Then, the induction hypothesis for the case of $n=2$ implies that the system in \eqref{eq:jordan_n_2} is ensemble controllable on $C(K,\mathbb{R})\oplus C(K,\mathbb{R}^k)=C(K,\mathbb{R}^{k+1})$ if and only if
$${\rm Im}(B(\b))={\rm Im}\left(\left[\begin{array}{c} b_1(\b) \\ B_k(\b) \end{array}\right]\right)={\rm Im}(b_1(\b))\oplus {\rm Im}(B_k(\b))=\mathbb{R}\oplus\mathbb{R}^k=\mathbb{R}^{k+1}$$ 
for all $\b\in K$, or equivalently, rank$(B(\b))=k+1$ for all $\b\in K$.
\hfill$\Box$

Note that, alternatively, Proposition \ref{prop:jordan_block} can also be proved using a Lie-algebraic approach, i.e., by showing that $\overline{\mathcal{L}}=C(K,\mathbb{R}^n)$ if and only if ${\rm rank}(B(\beta))=n$ for all $\b\in K$, where $\mathcal{L}$ is the Lie algebra generated by the drift and control vector fields of the system. 
This proposition also implies that a linear ensemble system with the system matrix similar to an $n$-by-$n$ Jordan block requires at least $n$ independent control inputs to be ensemble controllable on $C(K,\mathbb{R}^n)$. This result coincides with the case in which the system matrix is similar to a diagonal matrix with $n$ identical eigenvalue value functions (by Theorem \ref{thm:linear_diagonalizable_injective}). This observation reveals that coupling between different states in the Jordan block case does not necessarily promote 
controllability of linear ensemble systems.

On the contrary, in the classical case, a finite-dimensional linear system may be easier to be controllable if the state variables are coupled, which is illuminated by the following examples. Let's start with a system without state-coupling, i.e., the system matrix is diagonal with identical elements (eigenvalues), 
$$\frac{d}{dt}x(t)=Ax(t)+BU(t),$$
where $A={\rm diag}(\l,\dots,\l)\in\mathbb{R}^{n\times n}$. It is easy to verify that this system is controllable on $\mathbb{R}^n$ if and only if rank$(B)=n$, since the controllability Gramian $W=[B\mid \l B\mid \cdots\mid \l^{n-1}B]$ and $\text{rank}(W)=\text{rank}(B)$. In other words, this system must have at least $n$ controls so that each state receives an independent control to reach controllability. However, if the state variables are coupled, it is possible to reduce the number of control inputs required to guarantee controllability. In fact, one control signal is enough to warrant controllability of a linear system with the system matrix similar to a Jordan block, i.e., 
\begin{align}
	\label{eq:jordan_single}
	\frac{d}{dt}x(t)=Jx(t)+bu(t).
\end{align}
where 
$$J=\left[\begin{array}{ccccc} \l & 1 &  & & \\  & \l & 1 & & \\  &  & \ddots & \ddots &  \\  & & & \ddots & 1\\ &  &  & & \l \end{array}\right]\quad\text{ and }\quad b=\left[\begin{array}{c} 0 \\ 0 \\ \vdots \\ 1 \end{array}\right].$$
In this case, the controllability Gramian,
$$W=[b\mid Jb\mid \cdots\mid J^{n-1}b]=\left[\begin{array}{cccc} 0 & 0 & \cdots & 1 \\  \vdots & \vdots & \ddots & \vdots \\ 0 & 1 & \cdots & (n-1)\l^{n-2} \\ 1 & \l & \cdots & \l^{n-1} \end{array}\right],$$
is full rank.

\begin{figure}[t]
     \centering
     \includegraphics[width=0.9\columnwidth]{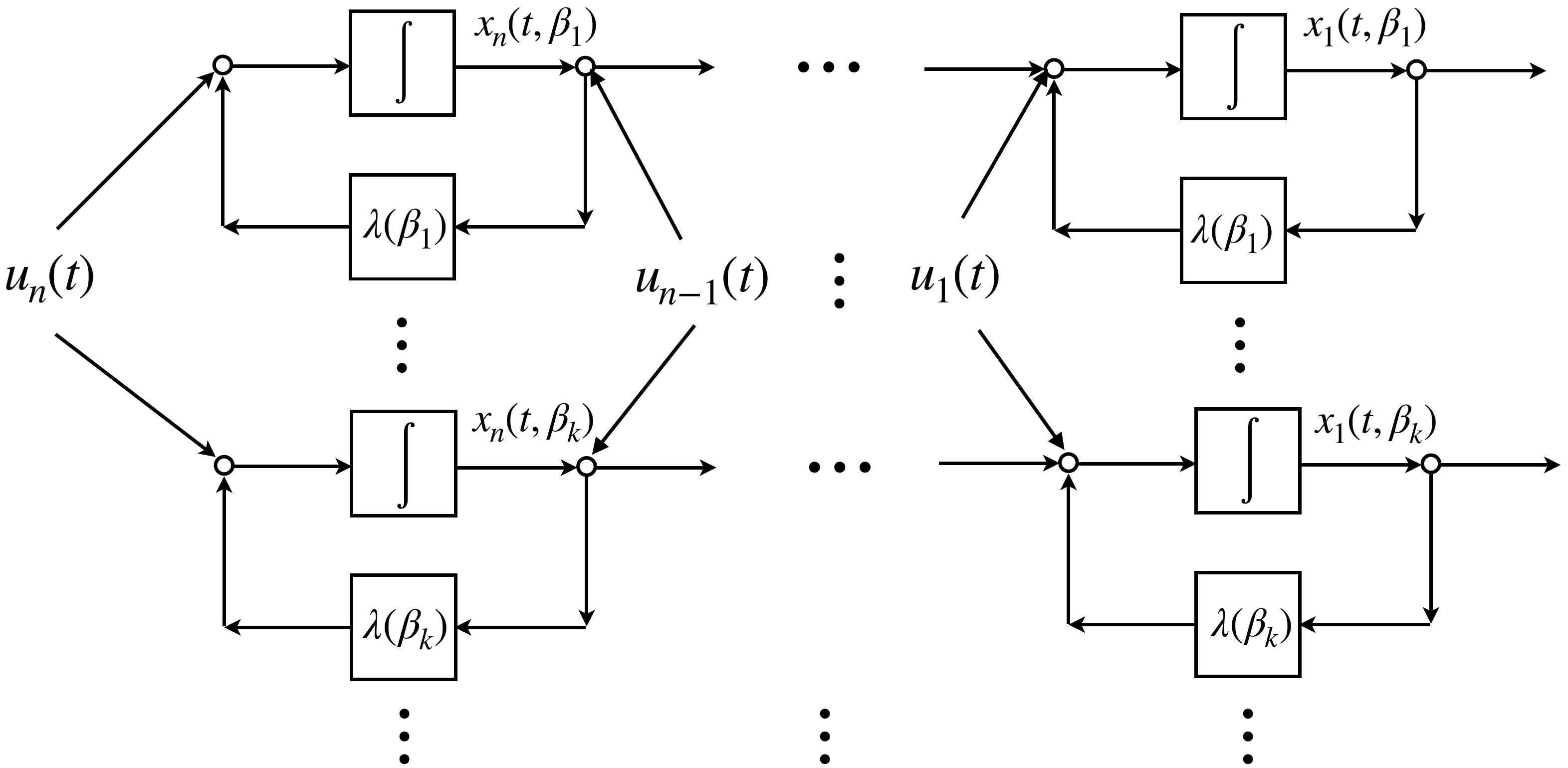}
     \caption{The illustration of an ensemble of $n$-dimensional linear systems with the drift in the form of a Jordan block using a population of chains of $n$ integrators.}
     \label{fig:integrators}
\end{figure}

\begin{remark}
	\label{rem:coupling}
	\rm
	The above example and the case studied in Proposition \ref{prop:jordan_block} illustrate the fundamental difference between how a control signal affects a classical linear and a linear ensemble system. For the classical control system as in \eqref{eq:jordan_single}, one input applied to a state variable (the last one) is sufficient to control the entire system, since the effect of this input is propagated to the other states through the coupling between them, while this is not the case for an ensemble system, which in fact requires at least $n$ independent 
	control fields to attain ensemble controllability. In practice, the system in \eqref{eq:jordan_single} represents the control of a chain of integrators, and the ensemble counterpart models the same control system in the presence of parameter uncertainty, i.e., $\b\in K$, or corresponds to the control of a collection of chains of $n$ integrators as shown in Figure \ref{fig:integrators}. For a single chain, 
	one control input applied to the last 
	integrator can be propagated to all other integrators sequentially so that the entire chain is controllable. However, for the ensemble case, 
	the control input robust to a state variable, e.g., $x_i(t,\b)$, will not be robust to the others, and hence having $n$ independent inputs are necessary and sufficient to fully control the ensemble chains as proved in Proposition \ref{prop:jordan_block}.
\end{remark}

The idea presented in Proposition \ref{prop:jordan_block} for analyzing controllability of ensemble systems with non-diagonalizable system matrix through the diagonalizable case lays the foundation for the study of more general ensemble systems.

\begin{corollary}
\label{cor:jordan}
	Given the time-invariant linear ensemble system, $\frac{d}{dt}X(t,\b)=A(\b)X(t,\b)+B(\b)U(t)$, 
	indexed by the parameter $\b$ taking values on a compact set $K\subset\mathbb{R}$, where $X(t,\cdot)\in C(K,\mathbb{R}^n)$, $U:[0,T]\rightarrow\mathbb{R}^m$ is piecewise constant, $A\in C(K,\mathbb{R}^{n\times n})$ is similar to a Jordan canonical form $J(\b)\in C(K,\mathbb{R}^{n\times n})$, and $B\in C(K,\mathbb{R}^{n\times m})$. Let $P^{-1}(\b)A(\b)P(\b)=\Lambda(\b)$ be the similar transformation of $A(t,\b)$ and let $Y(t,\b)=P(\b)X(t,\b)$ be the transformed system of $X(t,\b)$, obeying  
	$$\frac{d}{dt}Y(t,\b)=J(\b)Z(t,\b)+\widetilde{B}(\b)U(t),$$
	where $\widetilde{B}=P^{-1}B\in C(K,\mathbb{R}^{n\times m})$. Then, this system is uniformly ensemble controllable on $C(K,\mathbb{R}^n)$ if and only if the diagonalized system associated with $Y(t,\b)$, i.e., 
	$$\frac{d}{dt}Z(t,\b)=\Lambda(\b)Z(t,\b)+\widetilde{B}(\b)U(t),$$
	is ensemble controllable on $C(K,\mathbb{R}^n)$, where $\Lambda(\b)$ is the diagonal matrix whose elements are diagonal elements of $J(\b)$. 
\end{corollary}	
	{\it Proof.} This is a direct extension of Proposition \ref{prop:jordan_block}, and the proof follows that of Proposition \ref{prop:jordan_block} as well. 
	\hfill$\Box$

For a linear system, the off-diagonal entries in the system matrix represent the coupling between state variables. Although a Jordan matrix as $J$ in \eqref{eq:jordan_single} represents a special type of state-coupling, that is, each state is coupled only with its successor, the result presented in  Corollary \ref{cor:jordan} can be extended to the cases of general types of state-coupling depicted in the following theorem.

\begin{theorem}
\label{cor:summary}
	Given the time-invariant linear ensemble system, $\frac{d}{dt}X(t,\b)=A(\b)X(t,\b)+B(\b)U(t)$, 
	indexed by the parameter $\b$ taking values on a compact set $K\subset\mathbb{R}$, where $X(t,\cdot)\in C(K,\mathbb{R}^n)$, $U:[0,T]\rightarrow\mathbb{R}^m$ is piecewise constant, $A\in C(K,\mathbb{R}^{n\times n})$ has real eigenvalues with $\l_1(\b)$,$\dots,\l_n(\b)\in\mathbb{R}$ for all $\b\in K$, and $B\in C(K,\mathbb{R}^{n\times m})$. Let $$\frac{d}{dt}Y(t,\b)=T(\b)Y(t,\b)+\widetilde{B}(\b)U(t)$$
	be the system obtained by transforming $A(\b)$ to an upper triangular matrix $T(\b)$. Then, the following properties are equivalent: \\
	(i) This system is uniformly ensemble controllable on $C(K,\mathbb{R}^n)$.\\
	(ii) The diagonalized system 
		$$\frac{d}{dt}Z(t,\b)=\Lambda(\b)Z(t,\b)+\widetilde{B}(\b)U(t),$$
		where $\Lambda(\b)={\rm diag}(\l_1(\b),\dots,\l_n(\b))$, is ensemble controllable on $C(K,\mathbb{R}^n)$.\\
	(iii) The induced system of $Z(t,\b)$ parameterized by $\eta_1=\l_1(\b)$, $\dots$, $\eta_n=\l_n(\b)$, i.e., 
	$$\frac{d}{dt}\left[\begin{array}{c} Z_1(t,\eta_1) \\ \vdots \\ Z_n(t,\eta_n) \end{array}\right]=\left[\begin{array}{ccc} \eta_1I_{\kappa_1(\eta_1)} &  &   \\ & \ddots & \\ &   & \eta_nI_{\kappa_n(\eta_n)} \end{array}\right]\left[\begin{array}{c} Z_1(t,\eta_1) \\ \vdots \\ Z_n(t,\eta_n) \end{array}\right]+\left[\begin{array}{c} D_1(\eta_1) \\ \vdots \\ D_n(\eta_n) \end{array}\right]U(t),$$
	is controllable on $\mathbb{R}^N$ for each $n$-tuple $(\eta_1,\dots,\eta_n)\in K_1\times\cdots\times K_n$ with $K_i=\l_i(K)$, $i=1,\ldots,n$, where $\kappa_i(\eta_i)=|\l_i^{-1}(\eta_i)|$ is the cardinality of the inverse image of $\eta_i$ under $\lambda_i$, $N=\sum_{i=1}^n\kappa_i(\eta_i)$, $I_{\kappa_i(\eta_i)}$ is the $\kappa_i(\eta_i)\times\kappa_i(\eta_i)$ identity matrix, and $D_i(\eta_i)\in\mathbb{R}^{\kappa_i(\eta_i)\times m}$ is the ensemble controllability Gramian associated with the $i^{\rm th}$ state of the ensemble $Z(t,\b)$, i.e.,
	$$\frac{d}{dt}z_i(t,\b)=\l_i(\b)z_i(t,\b)+\tilde{b}_i(\b)U(t)$$
	with $\tilde{b}_i(\b)$ the $i^{\rm th}$ row of $\widetilde{B}(\beta)$.
\end{theorem}

{\it Proof.} The proof follows directly from Corollary \ref{cor:jordan}, Theorem \ref{thm:linear_ensemble}, and the observation that any triangular matrix shares the same eigenvalues as its Jordan canonical form.  
\hfill$\Box$

\begin{example}\rm
	Consider the linear ensemble system with the parameter $\b$ varying on $[0,1]$, given by
	\begin{align}
		\label{eq:coupled}
		\frac{d}{dt}X(t,\b)=\left[\begin{array}{cc} \b & 1 \\ 0 & \b^2\end{array}\right]X(t,\b)+\left[\begin{array}{ccc} 1 & 0 & 0 \\ 0 & 1 & \b \end{array}\right]U(t).
	\end{align} 
	We first observe that $A(\b)$ is diagonalizable when $\b\in (0,1)$ and is a Jordan block when $\b=0$ or $\b=1$. This indicates that there exists no similar transformation $P(\b)$ that will diagonalize $A(\b)$ or put it into a Jordan form for all $\b\in [0,1]$. Therefore, Theorem \ref{thm:linear_ensemble} and Proposition \ref{prop:jordan_block} can not be applied. However, because $A(\b)$ is upper triangular, according to Theorem \ref{cor:summary}, the system in \eqref{eq:coupled} is ensemble controllable on $C([0,1],\mathbb{R}^2)$ if and only if the diagonalized system,
	\begin{align}
	\label{eq:diagonal}
	\frac{d}{dt}X(t,\b)=\left[\begin{array}{cc} \b & 0 \\ 0 & \b^2\end{array}\right]X(t,\b)+\left[\begin{array}{ccc} 1 & 0 & 0 \\ 0 & 1 & \b \end{array}\right]U(t),
	\end{align}
	is ensemble controllable on $C([0,1],\mathbb{R}^2)$. By Theorem \ref{thm:linear_ensemble}, the system in \eqref{eq:diagonal} is uniformly ensemble controllable if its reparameterized system 
	\begin{align}
		\label{eq:single}
		\frac{d}{dt}Z(t,\eta_1,\eta_2)=\left[\begin{array}{cc} \eta_1 & 0 \\ 0 & \eta_2 I_{\kappa(\eta_2)} \end{array}\right]Z(t,\eta_1,\eta_2)+\left[\begin{array}{c} D_1(\eta_1) \\ D_2(\eta_2) \end{array}\right]U(t),
	\end{align}
	is controllable on $\mathbb{R}^N$ for each pair $(\eta_1,\eta_2)\in[0,1]\times[0,1]$, where $D_1(\eta_1)=\left[\begin{array}{ccc} 1 & 0 & 0 \end{array}\right]$ for all $0\leq\eta_1\leq1$,
	\begin{align*}
	D_2(\eta_2)=
	\begin{cases}
	\left[\begin{array}{ccc} 0 & 1 & 0 \end{array}\right],\ \qquad\eta_2=0, \\
	\left[\begin{array}{ccc} 0 & 1 & \sqrt{\eta_2} \\  0 & 1 & -\sqrt{\eta_2} \end{array}\right],\ 0<\eta_2\leq1,
	\end{cases}
	\text{and}\quad
	N=\kappa(\eta_2)+1=
	\begin{cases}
	2,\ \eta_2=0, \\
	3,\ 0<\eta_2\leq1.
	\end{cases}
	\end{align*}
	It is straightforward to check that the system in \eqref{eq:single} is indeed controllable for any choice of $0\leq\eta_1,\eta_2\leq1$, which implies uniform ensemble controllability of the system in \eqref{eq:coupled} on $C([0,1],\mathbb{R}^2)$.

\end{example}

It is worth noting that Theorem \ref{cor:summary} is a summary of the main results developed in this paper, which concludes the discussion about uniform ensemble controllability of the time-invariant linear ensemble system with parameters varying on a compact set and with the system matrix having real eigenvalues. The following closing remark of the paper then provides a vision of extending the scope our discussion to linear ensemble systems defined on broader classes of function spaces.

\begin{remark}[Amenable Extensions]
	\label{rmk:L_p}
	\rm
	According to the Stone-Weierstrass theorem (Theorem \ref{thm:Weierstrass} in Appendix \ref{appd:S-W_thm}), the ability to separate 
	points from different injective branches and from an overlapping spectrum guarantees that the Lie algebra generated by the drift and control vector fields is dense in $C(K,\mathbb{R}^n)$, which is the key connection of the notion of separating points to ensemble controllability.
	Moreover, because $C(K,\mathbb{R}^n)$ is dense in $L^p_n(K,\mu)$ for $1\leq p<\infty$ \cite{Folland99}, where $L^p_n(K,\mu)=\{f:K\rightarrow \mathbb{R}^n\mid \int_{K}\|f\|^pd\mu<\infty\}$, $\|\cdot\|$ is a norm on $\mathbb{R}^n$, and $\mu$ is a Radon measure on $K$, 
	the developed separating point methodology can also be directly utilized to analyze ensemble controllability of the time-invariant linear ensemble system defined on $L^p_n(K,\mu)$. Another possible extension of this work is to relax the compactness of the parameter space. Because the Weierstrass theorem also works for $C_0(\Omega,\mathbb{R}^n)$, the space of functions vanishing at infinity, in the case in which the parameter space $\Omega$ of an ensemble is not compact, the same results developed above in this paper are valid for ensemble systems 
	defined on $C_0(\Omega,\mathbb{R}^n)\subset C(\Omega,\mathbb{R}^n)$. 
\end{remark}

\section{Conclusion}
In this paper, we study the control of the time-invariant linear ensemble system defined on an infinite-dimensional space of continuous functions. We exploited the notion of separating points and applied techniques of polynomial approximation and Lie algebra to establish uniform ensemble controllability conditions. In particular, these conditions were constructed based on identifying the number of independent control fields sufficient to separate points in different injective branches and overlapping spectra of the drift dynamics, in which the analysis was carried out by the newly introduced ensemble controllability Gramian. We showed that ensemble controllability on an infinite-dimensional function space, e.g., $C(K,\mathbb{R}^n)$, can be examined through evaluating controllability of each individual system in the ensemble on a finite-dimensional space, e.g., $\mathbb{R}^n$. This development provides a unified and systematic procedure for analyzing controllability of control systems defined on an `infinite-dimensional' space by using a `finite-dimensional' method.

\appendix

\section{Stone-Weierstrass Theorem}
\label{appd:S-W_thm}
\begin{defn}[Separating Points] 
\label{def:separating_points}
	Let $X$ be a topological space. A subset $\mathcal{A}$ of $C(X,\mathbb{R})$ is said to separate points from $X$ if, for every $x,y\in X$ with $x\neq y$, there exists $f\in\mathcal{A}$ such that $f(x)\neq f(y)$ \cite{Folland99}. 
\end{defn}

Next, we state the Stone-Weierstrass Theorem using this terminology.

\begin{theorem}[Stone-Weierstrass Theorem] 
	\label{thm:Weierstrass}
	Let $X$ be a compact Hausdorff space. If $\mathcal{A}$ is a closed subalgebra of  $C(X,\mathbb{R})$ that separates points from $X$, then either $\mathcal{A}=C(X,\mathbb{R})$ or $\mathcal{A}=\{f\in C(X,\mathbb{R}):f(x_0)=0\}$ for some $x_0\in X$. The first alternative holds if and only if $\mathcal{A}$ contains the constant functions.
\end{theorem}

{\it Proof.} See \cite{Folland99}. \hfill$\Box$

By requiring $\mathcal{A}$ to be a subalgebra of $C(K,\mathbb{R})$, it means that $\mathcal{A}$ is a vector subspace of $C(K,\mathbb{R})$ such that $fg\in\mathcal{A}$ whenever $f,g\in\mathcal{A}$.  Furthermore, the Stone-Weierstrass theorem can be extended to the space $C_0(X,\mathbb{R})$ as the following corollary shows, where $X$ is not necessarily compact and $C_0(X,\mathbb{R})$ denotes the space of real-valued functions vanishing at infinity, i.e., the set of functions such that $\{x\in X:|f(x)|\geq\varepsilon\}$ is compact for any $f\in C_0(K,\mathbb{R})$ and $\varepsilon>0$.

\begin{corollary}[Stone-Weierstrass Theorem] 
	Let $X$ be a locally compact Hausdorff space. If $\mathcal{A}$ is a closed subalgebra of  $C_0(X,\mathbb{R})$ that separates points from $X$, then either $\mathcal{A}=C_0(X,\mathbb{R})$ or $\mathcal{A}=\{f\in C_0(X,\mathbb{R}):f(x_0)=0\}$ for some $x_0\in X$. 
\end{corollary}

{\it Proof.} See \cite{Folland99}. \hfill$\Box$

If $C_c(X,\mathbb{R})$ denotes the set of compactly supported real-valued functions defined on $X$, then it can be shown that $C_0(X,\mathbb{R})$ is the closure of $C_c(X,\mathbb{R})$ under the uniform topology.

\section{Construction of the Parameterization in Lemma \ref{lem:linear_2d}}
\label{app:2-d}

Recall the 2-dimensional linear ensemble system in \eqref{eq:linear_ensemble_2d}, i.e.,
$$\frac{d}{dt}\left[\begin{array}{c} y_1(t,\b) \\ y_2(t,\b) \end{array}\right]=\left[\begin{array}{cc} \l_1(\b) & 0 \\ 0 & \l_2(\b) \end{array}\right]\left[\begin{array}{c} y_1(t,\b) \\ y_2(t,\b) \end{array}\right]+\left[\begin{array}{c} \tilde{b}_1(\b) \\ \tilde{b}_2(\b) \end{array}\right]V(t),$$
where $\l_i\in C(K,\mathbb{R})$ and $\tilde{b}_i(\b)\in C(K,\mathbb{R}^m)$ for $i=1,2$, then we claim that there exists a parameterization $\psi: C(K',\mathbb{R})\rightarrow C(K,\mathbb{R}^2)$ such that the dynamics of  $z(t,\eta)\doteq\, \psi^{-1}\left[\begin{array}{c} y_1(t,\b) \\ y_2(t,\b) \end{array}\right]$ follows the one-dimensional linear ensemble system in \eqref{eq:1-d-param}, i.e.,
$$\frac{d}{dt}z(t,\eta)=\varphi(\eta)z(t,\eta)+b(\eta)V(t),$$
where $K'$ is a disjoint union of the two sets $K_1'$ and $K_2'$ that are homeomorphic to $K_1=\l_1(K)$ and $K_2=\l_2(K)$, respectively; $\varphi:K'\rightarrow K_1\cup K_2$ is continuous with $\varphi^{-1}(K_i)=K_i'$ for $i=1,2$ and satisfies 
$$
	\left\{\begin{array}{ll} |\varphi^{-1}(s)|=2, & s\in K_1\cap K_2\neq\varnothing, \\ |\varphi^{-1}(s)|=1, & s\in K_1\backslash K_2 \ \text{or} \ s\in K_2\backslash K_1;
	\end{array}\right.
$$
and the control vector field 
$b(\eta)=\tilde{b}_1(\l_1^{-1}\circ\varphi(\eta))$ for $\eta\in K_1'$ and $b(\eta)=\tilde{b}_2(\l_2^{-1}\circ\varphi(\eta))$ for $\eta\in K_2'$. This claim immediately follows the existence of $\varphi$ by defining $\psi: C(K',\mathbb{R})\rightarrow C(K,\mathbb{R}^2)$ as
$$f\mapsto\left[\begin{array}{c} f\circ(\varphi|_{K_1'})^{-1}\circ\l_1 \\ f\circ(\varphi|_{K_2'})^{-1}\circ\l_2 \end{array}\right].$$

If $K_1\cap K_2=\varnothing$, let $K_1'=K_1$ and $K_2'=K_2$, then one choice of $\varphi$ is the identity function on $K_1\cup K_2$. Next, we assume $K_1\cap K_2\neq\varnothing$. In this case, let $a=\sup(K_1\cap K_2)$ and $b=\inf(K_1\cap K_2)$, and define $K_1'=K_1-|a-b|=\{\eta-|a-b|:\eta\in K_1\}$ and $K_2'=K_2-|a-b|=\{\eta-|a-b|:\eta\in K_2\}$, then $K_1'$ and $K_2'$ are disjoint and homeomorphic to $K_1$ and $K_2$, respectively. Following that, let $id_{K_i}:K_i\rightarrow K_i$ denote the identity function on $K_i$ for $i=1,2$, then the function $\varphi:K_1'\cup K_2'\rightarrow K_1\cup K_2$ given by
$$\varphi(\eta)=\begin{cases} id_{K_1}(\eta+|a-b|),\quad\text{if }\eta\in K_1' \\ id_{K_2}(\eta-|a-b|),\quad\text{if }\eta\in K_2'\end{cases}$$ 
satisfies the requirements.

\bibliographystyle{siam}
\bibliography{SICON_SeparatingPt}

\end{document}